\documentclass[11pt]{article}
\usepackage[top=1.1in, left=1in, right=1in, bottom=1.1in]{geometry}
\usepackage{amsmath,amssymb,amsthm,amsfonts,verbatim}
\usepackage{enumerate}
\usepackage{microtype}
\usepackage{url}
\usepackage{algpseudocode}
\usepackage{tikz}
\usetikzlibrary{arrows}
\usepackage{amsmath,amssymb}
\usepackage[justification=centering]{caption}
\usepackage[justification=centering]{subcaption}
\usepackage[english]{babel}
\usepackage{array}
\usepackage{psfrag}
\usepackage{graphicx}
\usepackage{color}
\usepackage{setspace}
\usepackage{hyperref}
\usepackage{algorithm2e}
\usepackage[utf8]{inputenc}
\usepackage{cleveref}
\usepackage{mathrsfs} % script font
\usepackage{authblk}
\usepackage{threeparttable}
\usepackage[titletoc,title]{appendix}
\usepackage{xcolor}
\usepackage{xfrac}
\usepackage{caption}
\SetKw{KwInitialization}{Initialization}
\SetKw{KwTheMainBody}{The main body}

\newtheorem{lemma}{Lemma}[section]

\newtheorem{assumption}{Assumption}
\newtheorem{theorem}{Theorem}[section]

\newtheorem{definition}[theorem]{\bf Definition}

\newtheorem{remark}{Remark}[section]

\definecolor{red}{rgb}{1,0.2,0.2}

\def \mr {\mathbb{R}}

\def \mx {\mathbf{x}}
\def \ml {\mathbf{l}}
\def \mw {\mathbf{w}}
\def \mb {\mathbf{b}}
\def \my {\mathbf{y}}
\def \mz {\mathbf{z}}

\def \myp {\mathbf{p}}

\SetKw{KwInitialization}{Initialization}
\SetKw{KwTheMainBody}{The main body}

\hypersetup{colorlinks=false,linkbordercolor=red,linkcolor=blue,pdfborderstyle={/S/U/W 1}}

\makeatletter
\newcommand*{\rom}[1]{\expandafter\@slowromancap\romannumeral #1@} % Roman numbers
\makeatother

\hypersetup{colorlinks=true,linkbordercolor=red,linkcolor=blue,pdfborderstyle={/S/U/W 1}}

\title{Computing the quasipotential for highly dissipative and chaotic SDEs. An application to stochastic Lorenz'63.}%
\author[1]{Maria Cameron\thanks{cameron@math.umd.edu}}
\author[1]{Shuo Yang\thanks{shuoyang@math.umd.edu}}
\affil[1]{Department of Mathematics, University of Maryland, College Park, MD 20742, USA}
 
\begin{document}

\maketitle

%\keywords{
%quasipotential,
%ordered line integral method,
%Lorenz'63,
%maximum likelihood transition path.
%}

%\subjclass[2010]{65N99, 65P99, 58J65}

\begin{abstract}
The study of noise-driven transitions occurring rarely on the time-scale of systems modeled by SDEs is of crucial importance for understanding such phenomena as genetic switches in living organisms and magnetization switches of the Earth. For a gradient SDE, the predictions for transition times and paths between its metastable states are done using the potential function. For a nongradient SDE, one needs to decompose its forcing into a gradient of the so-called quasipotential and a rotational component, which cannot be done analytically in general.

We propose a methodology for computing the quasipotential for highly dissipative and chaotic 
systems built on the example of Lorenz'63 with an added stochastic term.
It is based on the ordered line integral method, a Dijkstra-like quasipotential solver, and
combines 3D computations in whole regions, a dimensional reduction technique, and 2D computations on radial meshes on manifolds or their unions. Our collection of source codes is available on M. Cameron's web page and on GitHub.
\end{abstract}

%%%%%%%%%%%%%%%%%%%%%%%%%%%%%%%%%%%%%%%%%%

\section{Introduction}
Suppose a system is evolving according to a stochastic differential equation (SDE) of the form 
\begin{equation}
\label{sde0}
d\mx = \mb(\mx)dt + \sqrt{\epsilon}d\mathbf{w}, \quad \mx\in \mr^d,
\end{equation}
where $\mb(\mx)$ is a continuously differentiable vector field, $d\mathbf{w}$ is the standard Brownian motion, and $\epsilon$ is a small parameter.
The quasipotential is a key function of the large deviation theory (LDT) \cite{FW} that allows one to find a collection of useful asymptotic estimates
for long-time dynamics of such systems. They include
 the invariant probability measure, expected escape times from neighborhoods of 
 attractors of the corresponding ODE $\dot{\mx} = \mb(\mx)$
 lying within their basins, 
and maximum likelihood escape paths from the basins. 
The quasipotential can be viewed as an analogue to the 
potential function $V(\mx)$, $\mx\in \mr^d$, for a gradient SDE with deterministic term $-\nabla V(\mx)$.
The quasipotential is defined as the solution to the Freidlin-Wentzell action  functional minimization problem. 
The quasipotential is Lipschitz-continuous in any bounded domain but not necessarily continuously differentiable \cite{quasi}.
Unfortunately, it can be found analytically only in special cases, for example, for linear SDEs \cite{chen,chen1}.

Ordered line integral methods (OLIMs) for computing the quasipotential for  SDEs of the form \eqref{sde0}
 in whole regions on regular rectangular meshes were introduced in \cite{DM} for 2D and extended to 3D in \cite{YPC}. 
 They are Dijkstra-like solvers that advance the solution from mesh points with smaller values to those with larger values\footnotemark[1]
 without iteration. 
 \footnotetext[1]{This is only approximately true. See Ref. \cite{OUM2003} for details.}
Their general structure is inherited from the ordered upwind method (OUM) \cite{OUM2001,OUM2003}, but there are important differences.
First, unlike the OUM that uses the upwind finite difference scheme, the OLIMs solve a local functional minimization problem at every step
approximating a segment of curve with a segment of straight line, 
and the integral along it by an at least second order accurate quadrature rule. 
This renders their observed rate of convergence superlinear for some cases, and reduces error constants 
by two to three orders of magnitude in comparison with the OUM.
Second, while the OUM is practical only for 2D problems due to large CPU times in larger dimensions, 
the OLIMs have been successfully extended for 3D. 
This became possible due to the hierarchical update strategy \cite{DM,YPC}, the use of the Karush-Kuhn-Tucker optimality 
conditions to eliminate unnecessary updates,
and a number of implementational rationalizations. 

In previous works  \cite{DM,DM1,YPC},  the OLIMs were developed for computing the quasipotential for mild-to-moderate ratio 
$\Xi(\mx)$ of the magnitudes of the rotational and potential components of the vector field $\mb(\mx)$ in \eqref{sde0}. 
In all test problems considered  in \cite{DM,DM1,YPC}, $\Xi(\mx)$ did not exceed 10 within 
in the important region around the attractor with respect to which the quasipotential was computed.
For all these test problems, the black-box algorithms \cite{DM,DM1,YPC} produced numerical solutions with small relative errors.

%A salient parameter of the OLIMs is the update factor $K$ that defines the radius $Kh$ ($h$ is the mesh step)
%of neighborhood of any mesh point
%from which it can be updated. $K$ is an integer that depends on the mesh size and the ratio $\Xi(\mx)$ of the magnitudes of the 
%rotational and potential components of the vector field $\mb(\mx)$ in \eqref{sde0}. This ratio is not known in 
%advance but can be estimated in some cases and can be found from an accurate enough computed solution.
%In all test problems and applications considered in \cite{DM,DM1,YPC}, $\Xi(\mx)$ did not exceed 10 within 
%in the important region around the attractor with respect to which the quasipotential was computed.
%%,the level sets of the quasipotential completely lying within the basin of the attractor with respect to which the quasipotential was computed,
%%perhaps except for those level sets that include some small neighborhood of the transition 
%%state in the case where the transition state is the saddle limit cycle.
%For all these test problems, the black-box algorithms \cite{DM,DM1,YPC} produced numerical solutions with small relative errors.

Unfortunately, if one applies the black-box {\tt olim3D} quasipotential solver from \cite{YPC} 
to a highly dissipative and chaotic system such as Lorenz'63 with an added small white noise, the 
 relative error of the numerical solution might be large leading to completely wrong estimates for escape rates.
For the parameter values $\sigma = 10$, $\beta = \sfrac{8}{3}$, and $\rho\gtrsim 15$, the quasipotential computed 
with respect to one of the point attractors
 will become progressively inaccurate as $\rho$ increases. We show in this work that, as $\rho$ approaches 
 $\rho_2\approx 24.74$ (where a subcritical Hopf bifurcation happens), 
 the upper bound for the ratio $\Xi(\mx)$ blows up at any point of the computational domain of interest. 
Even if one uses a very good desktop computer\footnotemark[2],
this problem  cannot be cured by mesh refinement due to the limited computer's memory: 
the size of a 3D mesh cannot exceed $1000^3$ by much.
\footnotetext[2]{We use iMac 2017 with processor 4.2 GHz Intel Core i7 and memory 64 GB 2400 MHz DDR4.}

In this work, we propose an approach for computing the quasipotential, finding maximum likelihood transition paths, and estimating 
escape times from basins of attractors for highly dissipative and possibly chaotic systems perturbed by small white noise.
This approach is suitable for systems where the 3D dynamics, 
after some short transition time, takes place in a small neighborhood of a 2D manifold  or a union of 2D manifolds
consisting of
certain characteristics of the corresponding ODE  (see Assumption \ref{assumption1} in Section \ref{sec:challenge} below).
Whether or not this phenomenon takes place can be identified from the plots of the 3D level sets of the  computed quasipotential.
We develop a technique for extracting these manifolds and generating so-called radial meshes on them. 
We adjust and test the OLIM for 2D radial meshes and compute the quasipotential on the constructed 2D manifolds or their unions.

The proposed techniques have been developed on the stochastic Lorenz'63:
\begin{equation}
\label{sde1}
d\mx = \left[\begin{array}{c}\sigma(x_2 - x_1)\\
x_1(\rho - x_3) - x_2\\
x_1x_2 - \beta x_3
\end{array}\right]dt + \sqrt{\epsilon}d\mathbf{w},\quad{\rm where}\quad \mx\equiv\left[\begin{array}{c}x_1\\x_2\\x_3\end{array}\right],
\end{equation}
with $\sigma = 10$, $\beta = \sfrac{8}{3}$,  and $0.5 \le \rho < \rho_2\approx 24.74$.
To the best of our knowledge, this is the first time when the quasipotential is computed 
for a chaotic 3D system in the whole region and 3D computations are refined by 2D computations on certain manifolds. 
We study transitions between the stable equilibria at $\rho = 12$, 15, and 20, and between the stable equilibria and the strange attractor at $\rho = 24.4$,
and find a collection of quasipotential barriers for them. 
Our transition paths obtained by a direct integration using the computed quasipotential
can be compared to those 
found in \cite{zhou} using the minimum action method, a path-based method consisting in a direct minimization of 
the Freidlin-Wentzell action in the path-space. 
At $\rho = 24.4$, we compare two plausible transition mechanisms from the strange attractor to the equilibria.
We offer a number of plots of 3D level sets of the quasipotential at various values 
of $\rho$ varying from $0.5$ to $24.4$ and supplement them with links to youtube videos for a better 3D visualization.
For $\rho\ge15$, when 2D approximation becomes accurate enough, we perform refined 2D computations of the quasipotential.

Aiming at making our results readily reproducible, we 
made most of the codes developed for this project
publicly available at M. Cameron's web page \cite{mariakc} -- see the package {\tt Qpot4lorenz63.zip}, and on GitHub \cite{github}.
All codes mentioned throughout this paper are included in this package.
A user guide for the codes is also provided there.  

The techniques developed in this work can be used for analysis of other stochastic systems.
For example, the computation of the quasipotential for the 3D genetic switch model from \cite{lv} would benefit from
performing a refined 2D computation on a radial mesh on a 2D manifold as suggested by Fig. 9 in \cite{YPC}.
Gissinger's 3D model \cite{giss} relevant for the reversals of the magnetic 
field of the Earth can be analyzed using the tools developed in this work. 

The rest of the paper is organized as follows. 
In Section \ref{sec:background}, some necessary background on the quasipotential is given.
A brief overview of the dynamics of Lorenz'63 at $\sigma=10$, $\beta = \sfrac{8}{3}$, and $0<\rho<\infty$
 is offered in Section \ref{sec:lorenz} and Appendix \ref{appLorenz}.
Numerical techniques for computing the quasipotential are described in Section \ref{sec:numerics}.
The application  to stochastic Lorenz'63 
is presented in Section \ref{sec:results}. 
We summarize our findings in Section \ref{sec:conclusions}. 
Some technical details are explained in Appendices \ref{sec:newapp}--\ref{sec:appB}.

%%%%%%%%%%%%%%%%%%%%%

\section{Definition and significance of the  quasipotential}
\label{sec:background}
To explain what is the quasipotential \cite{FW},  we first assume that the vector field $\mb(\mx)$  in SDE \eqref{sde0}
admits the following smooth orthogonal decomposition:
\begin{equation}
\label{b}
\mb(\mx) = -\frac{1}{2}\nabla u(\mx) +\ml(\mx),\quad  \nabla u(\mx)\cdot \ml(\mx) = 0.
\end{equation}
If  $\ml(\mx)\equiv \mathbf{0}$, i.e., if the field $\mb(\mx)$ were gradient,
the Gibbs measure 
\begin{equation}
\label{gibbs}
\mu(\mx) = Z^{-1}e^{-u(\mx)/\epsilon}
\end{equation} 
would be the invariant probability density for SDE \eqref{sde0}. Suppose $\ml(\mx)$ is not identically zero.
Plugging the Gibbs measure \eqref{gibbs} into the  stationary Fokker-Planck  equation for SDE \eqref{sde0}
\begin{equation}
\label{FP}
\frac{1}{2}\Delta \mu(\mx) -\nabla\cdot(\mu(\mx)\mb(\mx)) = 0
\end{equation}
we find that it is invariant if and only if $\ml(\mx)$ is divergence-free, i.e., $\nabla\cdot\ml(\mx)\equiv 0$.
In this case, the function $u(\mx)$ would play the role of a potential. 

Unfortunately, the orthogonal decomposition \eqref{b} where $\ml(\mx)$ is divergence-free does not typically exist. 
However, a function $U(\mx)$ called the quasipotential that gives asymptotic estimates for the invariant probability measure near attractors of 
$\dot{\mx}=\mb(\mx)$ in the limit $\epsilon\rightarrow 0$ can be designed \cite{FW}. 

Suppose that the vector field $\mb(\mx)$ is continuously differentiable. In addition, we assume that 
the ODE
\begin{equation}
\label{ode}
\dot{\mx}=\mb(\mx)
\end{equation}
has a finite number of attractors, and every trajectory of \eqref{ode} remains in a bounded region as $t\rightarrow \infty$.
Let $A$ be an attractor of \eqref{ode}.
The quasipotential with respect to  $A$ is defined as the solution of the minimization problem 
\begin{equation}
\label{qpot}
U(\mx) = \inf_{\phi,T_0,T_1}\left\{ S_{T_0,T_1}(\phi)~|~\phi(T_0)\in A,~\phi(T_1) = \mx\right\},
\end{equation}
 where the infimum of the Freidlin-Wentzell action 
 \begin{equation}
 \label{FWA}
S_{T_0,T_1}(\phi) = \frac{1}{2}\int_{T_0}^{T_1}\|\dot{\phi} - \mb(\phi)\|^2dt
\end{equation}
is taken over the set of absolutely continuous paths $\phi$ with endpoints at $A$ and $\mx$, and all times $T_0,T_1\in\mr$.
The infimum with respect to $T_0$ and $T_1$ can be taken analytically \cite{FW,hey1,hey2}
resulting in the geometric action (see Appendix \ref{sec:newapp})
\begin{equation}
\label{GA}
S(\psi) = \int_0^L\left(\|\psi'\|\|\mb(\psi)\|-\psi'\cdot\mb(\psi)\right)ds,
\end{equation}
where the path $\psi$ is parametrized by its arclength, and $L$ is the length of $\psi$.
As a result, the definition of the quasipotential can be rewritten in terms of the geometric action:
\begin{equation}
\label{qpot1}
U(\mx) = \inf_{\psi}\left\{ S(\psi)~|~\psi(0)\in A,~\psi(L) = \mx\right\}.
\end{equation}
We have been using definition \eqref{qpot1} to develop quasipotential solvers.

Using Bellman's principle of optimality \cite{bellman}, one can show \cite{quasi} that the quasipotential $U(\mx)$ satisfies the Hamilton-Jacobi equation (see Appendix \ref{sec:newapp})
\begin{equation}
\label{HJ}
\frac{1}{2}\|\nabla U(\mx)\|^2 + \mb(\mx)\cdot\nabla U(\mx) = 0, \quad U(A) = 0.
\end{equation}
Eq. \eqref{HJ} implies that 
\begin{equation}
\label{ort}
\mb(\mx) = -\frac{1}{2}\nabla U(\mx) + \ml(\mx),~~{\rm where}~~\ml(\mx): = \mb(\mx) + \frac{1}{2}\nabla U(\mx)~\text{is orthogonal to}~ \nabla U(\mx).
\end{equation}
We will refer to $-(\sfrac{1}{2})\nabla U(\mx) $ and $\ml(\mx)$ as the potential and rotational components respectively. 

We remark that the boundary value problem (BVP) \eqref{HJ} is ill-posed.
It always has the trivial solution identically equal to zero
and may or may not have a smooth nontrivial solution. 
The quasipotential defined by \eqref{qpot} or \eqref{qpot1}
is a viscosity solution\footnotemark[1]  to \eqref{HJ} \cite{crandall}.
\footnotetext[1]{A viscosity solution to a first-order nonlinear PDE $f(\mx,u,\nabla u) =0$ is a continuous but possibly nondifferentiable function
obtained as the limit of a sequence of smooth solutions to $f(\mx,u,\nabla u) =\epsilon\Delta u$ as $\epsilon\rightarrow\infty$.}
The other complication is that even a nontrivial solution to this BVP, classical or viscosity, may not be unique
due to the fact that the boundary condition is imposed on an attractor \cite{ishii}. 
For example, if $\mb(\mx) = B\mx$  where $B$ is a matrix with all eigenvalues having negative real parts,  
the number of solutions of \eqref{HJ} with the BC $u(\mathbf{0}) = 0$ is equal to the number of invariant subspaces for $B$.

Nonetheless, \eqref{HJ} is instrumental in deriving the equation 
for \emph{minimum action paths} (MAPs) a.k.a. \emph{maximum likelihood paths}  or \emph{instantons}
that minimize the geometric action \eqref{GA} \cite{FW,quasi} (see Appendix \ref{sec:newapp}):
\begin{equation}
\label{psi}
\psi'(s) = \frac{\mb(\psi(s)) + \nabla U(\psi(s))}{\|\mb(\psi(s)) + \nabla U(\psi(s))\|}.
\end{equation}
Once the quasipotential is computed, 
one can shoot a MAP from a given point $\mx$ back to the attractor $A$ by integrating \eqref{psi} backward in $s$.
Alternatively, MAPs can be found by path-based methods \cite{mam,zhou1,hey1,hey2}
that directly minimize the Freidlin-Wentzell action or the geometric action.

The mentioned asymptotic estimate for the invariant probability density within a level  set of the quasipotential completely lying in the 
basin $\mathcal{B}(A)$ of $A$ is \cite{FW}
\begin{equation}
\label{amu}
\mu(\mx)\asymp e^{-U(\mx)/\epsilon},\qquad~~{\rm i.e.}\qquad~~\lim_{\epsilon\rightarrow 0}\left(-\epsilon\log\mu(\mx)\right) = U(\mx).
\end{equation}
The symbol $\asymp$ denotes the logarithmic equivalence clarified in \eqref{amu}.
The expected escape time from $\mathcal{B}(A)$ can also be estimated up to exponential order \cite{FW}:
\begin{equation}
\label{etime}
\mathbb{E}[\tau_{\mathcal{B}(A)}] \asymp e^{U(\mx^{\ast})/\epsilon},~~{\rm where}~~U(\mx^{\ast}) = \min_{\mx\in\partial \mathcal{B}(A)}U(\mx).
\end{equation}
In some common special cases, a sharp estimate for the expected escape time can be obtained \cite{bouchet}.

The term \emph{transition state} is often encountered in  chemical physics literature. Mostly it refers to a saddle  
lying on the manifold separating two basins of attraction. 
The dynamics of the Lorenz system are complicated, and basins of its attractors are tightly interlaced for $\rho\gtrsim 20$.
To accommodate such situations, 
we will define the term \emph{escape state}.
\begin{definition}
Consider a system evolving according to SDE \eqref{sde0}. 
Let $A$ be an attractor of the corresponding ODE \eqref{ode}.
The escape state from $A$ is the set of points minimizing the quasipotential with respect to $A$ over the boundary of the basin of $A$.
\end{definition}
The quasipotential  at the escape state of $A$ defines the expected escape 
time from the basin of $A$ up to exponential order according to Eq.
\eqref{etime}.

%
%In summary, the quasipotential allows us to obtain important asymptotic estimates quantifying the dynamics of \eqref{g}.
%The quasipotential is rarely available analytically, however, 
%numerical methods have been developed for computing it in 2D and 3D \cite{DM,YPC}.
%In the next section, we will demonstrate how one can visualize the 
%dynamics of the stochastic Lorenz'63 system \eqref{sde1}  using the quasipotential.

%%%%%%%%%%%%%%%%%%%%%%%%%

\section{A brief overview of Lorenz'63 }
\label{sec:lorenz}
The Lorenz'63 system 
\begin{align}
\dot{x_1} &= \sigma(x_2 - x_1),\notag\\
\dot{x_2} & = x_1(\rho - x_3) - x_2,\label{lorenz}\\
\dot{x_3} & = x_1x_2 - \beta x_3 \notag
\end{align}
is one of the most fascinating  and transformative ODE models proposed in the twentieth century.
E. Lorenz \cite{lorenz} derived it from Saltzman's 2D cellular convection model \cite{saltzman} 
using a Fourier expansion and truncating the trigonometric series to include 
a total of three terms. He proved that the resulting system exhibits a new type of long-term behavior. 
All trajectories of \eqref{lorenz} stay in a bounded region.
For $\sigma = 10$, $\beta = \sfrac{8}{3}$, and $\rho = 28$, their $\omega$-limit sets form an ``infinite complex of surfaces", i.e., a fractal, 
whose Hausdorff dimension is 2.06   \cite{hdim}, 
later named the Lorenz attractor.  
The Lorenz map \cite{lorenz}, a 1D map  $z_{n+1}=f(z_n)$, where $z_n$ is the $n$th maximum of the $z$-component  of a trajectory, and $f$ is the function estimated numerically, 
explained the divergence of arbitrarily close characteristics.
It has become instrumental for analysis of chaotic dynamical systems.

The study of the Lorenz'63 system bursted in mid-1970s, perhaps due to the progress in the computer industry.
A number of remarkable properties and quantitative characteristics have been discovered.
The topological structure of the Lorenz attractor was studied in \cite{GW,Rand,Williams}. 
The phenomenon called preturbulence was described in \cite{preturbulence}. The value $\rho_1\approx 24.06$ at which the Lorenz attractor is born  
for $\sigma = 10$ and $\beta = \sfrac{8}{3}$ was found in \cite{yorke} using a functional fit to the Lorenz map. 
Homoclinic explosions, period-doubling cascades, and periodicity windows  were investigated in \cite{sparrow}.
A beautiful overview of the Lorenz system is given in \cite{strogatz} (Chapters 9--12).
Nowadays, the Lorenz system is a popular test model for new methods in such fields as 
machine learning and forecasting (e.g. \cite{dudul,sorrentino&ott,HBS}).

It is easy to check that \eqref{lorenz} is invariant under the symmetry transformation  $(x_1,x_2,x_3)\mapsto (-x_1,-x_2,x_3)$.
We fix the parameters $\sigma = 10$ and $\beta=\sfrac{8}{3}$ and consider the dynamics of \eqref{lorenz} as $\rho$ grows
from zero to infinity. The notation and bifurcations important for the further presentation are summarized in Table \ref{table:lorenz}.
A more detailed description of the dynamics of \eqref{lorenz} for $0<\rho<\infty$ is given in Appendix \ref{appLorenz}.

\begin{table}[htp]
\caption{A summary of bifurcations and notation for Lorenz'63 \eqref{lorenz} for $\sigma = 10$, $\beta = \sfrac{8}{3}$,
and $0<\rho\le \rho_2\approx 24.74$.
}
\begin{center}
\begin{tabular}{|c|c|}
\hline
{\bf Range of $\rho$}&{\bf Comments \& notation}\\
\hline
$0<\rho<1$ &  The origin is the unique globally attracting equilibrium.\\
\hline
$\rho = 1$ & Supercritical pitchfork bifurcation. \\
\hline
$1 < \rho < \rho_0\approx 13.926$ & The origin is a Morse index one saddle for $1<\rho<\infty$.  \\
& Equilibria $C_{\pm}$ are located at\\
&$
C_{\pm} = \left(\pm\sqrt{\beta(\rho-1)},\pm\sqrt{\beta(\rho-1)},\rho-1\right).
$
\\
& $C_{\pm}$ are asymptotically stable for $1<\rho<\rho_2$.\\
\hline
$\rho = \rho_0\approx 13.926$ & Homoclinic orbits starting and ending at the origin exist. \\
\hline
$\rho_0 < \rho <\rho_1\approx 24.06$ & $C_{\pm}$ are surrounded by saddle cycles $\gamma_{\pm}$ respectively.\\
&Chaotic dynamics (``preturbulence") is developing as $\rho$ grows. \\
& We introduce cones $\Upsilon_{\pm}$ with vertices at $C_{\pm}$ and \\
& passing through $\gamma_{\pm}$ respectively:\\
&$
\Upsilon_{+}: = \{ C_{+} + t(\mx - C_{+}) ~|~t\ge 0,~\mx\in\gamma_{+} \}.
$
 \\
\hline
$\rho = \rho_1\approx 24.06$ & The birth of the Lorenz attractor $A_L$ (a strange attractor).\\
\hline 
$\rho_1 < \rho <\rho_2\approx 24.74$ & $A_L$ coexists with asymptotically stable equilibria $C_{\pm}$. \\
\hline
$\rho = \rho_2\approx 24.74$ & A subcritical Hopf bifurcation: $\gamma_{\pm}$ shrink to $C_{\pm}$ respectively.\\
\hline
\end{tabular}
\end{center}
\label{table:lorenz}
\end{table}%

In this work, we consider the Lorenz system perturbed by small white noise \eqref{sde1}.
The noise term regularizes the 
chaotic deterministic dynamics of \eqref{lorenz} in the sense that one can predict the future probability 
density function given the current one by solving the Fokker-Planck equation.
On the other hand, the presence of the noise term enables escapes from any neighborhood of an attractor of \eqref{lorenz}.
If $\rho$ is such that there are multiple attractors, noise-induced transitions between their neighborhoods become possible.

%%%%%%%%%%%%%%%%%%%%%%%%%

\section{Numerical methods}
\label{sec:numerics}
In this section, we describe numerical techniques developed for computing the quasipotential for highly dissipative and chaotic systems where the ratio 
of  the magnitudes of the rotational and potential components is of the order of $10^3$.

%%%
\subsection{A brief overview of ordered line integral methods (OLIMs)}
We start with a brief overview the OLIMs. 
A comprehensive description of the implementation of the OLIM in 3D is provided in \cite{YPC}. 
It involves many technical details that are important for making the solver fast. 
A C source code {\tt olim3D4Lorenz63.c} set up to compute the quasipotential for \eqref{sde1}
and an instruction on how to run are available in \cite{mariakc,github}.

The OLIMs belong to the family of label-setting algorithms \cite{CV}
and inherit their set of labels from the OUM \cite{OUM2001,OUM2003}. 
Labels of mesh points indicate their statuses.
A mesh point is {\sf Accepted} if the value of the computed function (the quasipotential in our case) 
is finalized at it and all its nearest neighbors also have finalized values. 
{\sf Accepted } points are not used for updating values at other mesh points.
A mesh  point is {\sf Accepted Front} if the value at it is finalized but it has at least one nearest neighbor 
with a non-finalized value. {\sf Considered} mesh points are those with non-finalized tentative values that have at least one {\sf Accepted Front} nearest neighbor.
{\sf Unknown} mesh points have no {\sf Accepted Front} nearest neighbors and the values at them have not been proposed yet.

The OLIMs use  several kinds of neighborhoods of mesh points. 
The neighborhoods are defined via distances between indices of the mesh points. 
Let $\myp:=(i,j,k)\in\mathbb{Z}^3$ and $\myp_0:=(i_0,j_0,k_0)\in\mathbb{Z}^3$ be the lattice points corresponding to the mesh points 
$\mx$ and $\mx_0$ respectively. In other words, $\myp$ and $\myp_0$ are the indices of the mesh points $\mx$ and $\mx_0$, respectively.
Recall that the $l_q$, $q = 1,2$, and $l_{\infty}$ distances between $\myp$ and $\myp_0$ are defined as
\begin{align*}
\|\myp-\myp_0\|_{q} &: = \left[|i-i_0|^q + |j - j_0|^q + |k - k_0|^q\right]^{1/q}~~{\rm and}\\
\|\myp - \myp_0\|_{\infty} &: = \max\{ |i-i_0|, |j - j_0| ,|k - k_0|\},
\end{align*}
respectively. Let $\mathcal{I}$ be the set of indices of all mesh points.
\begin{itemize}
 \item
  The near neighborhood typically containing 26 points
 $$
 \mathcal{N}_{{\rm near}}(\myp_0): = 
  \{ \myp \in \mathcal{I}~|~\|\myp- \myp_0\|_1 \le 3 ~{\rm and}~\|\myp - \myp_0\|_{\infty} = 1\}
$$ 
is used for recruiting {\sf Unknown} points to {\sf Considered} and changing the status of {\sf Accepted Front} points to {\sf Accepted}. 
Correspondingly, the near neighborhood of the mesh point $\mx_0$ is defined as
$$
 \mathcal{N}_{{\rm near}}(\mx_0): = \left\{ \mx~|~ \myp\in \mathcal{N}_{{\rm near}}(\myp_0)\right\}.
$$ 
\item 
The far neighborhood  $\mathcal{N}_{{\rm far}}^K(\myp_0)$, where $K$ is the update factor (a positive integer chosen by the user), 
consists approximately\footnotemark[2] of all
lattice points $\myp\in\mathcal{I}$ such that $\myp\neq \myp_0$ and the $l_2$ distance $\|\myp-\myp_0\|_2\le K$.
It is used for updating {\sf Considered} points.
\footnotetext[2]{
More precisely, $\myp\in \mathcal{N}_{{\rm far}}^K(\myp_0)$ iff $\myp\neq \myp_0$, $\myp\in \mathcal{I}$, and
 $ |i - i_0|  \le K$,
$ |j - j_0|  \le {\rm ceil}\left(\sqrt{K^2 - |i - i_0|^2}\right)$, and
$ |k - k_0|  \le {\rm ceil}\left(\sqrt{K^2 - \min\{|i - i_0|^2 +|j - j_0|^2,K^2\} }\right)$.
Defined so, $\myp\in \mathcal{N}_{{\rm far}}^K(\myp_0)$ is slightly larger than $\{\myp\in\mathcal{I}~|~\myp\neq \myp_0,~\|\myp-\myp_0\|_2\le K\}$.
}
Correspondingly, the far neighborhood of the mesh point $\mx_0$ is defined as
$$
 \mathcal{N}^K_{{\rm far}}(\mx_0): = \left\{ \mx~|~ \myp\in \mathcal{N}^K_{{\rm far}}(\myp_0)\right\}.
$$ 
If the mesh steps in $x_i$, $i=1,2,3$, are all equal to $h$ then the far neighborhood of $\mx_0$ is approximately the ball centered
at $\mx_0$ of radius $Kh$.
\end{itemize}

At the start, all mesh points are {\sf Unknown}. Initialization consists in computing tentative 
values at the mesh points lying near the attractor, switching their status to {\sf Considered}, 
and adding them to the binary tree. The binary tree maintains the heap sort of the values at {\sf Considered} points so that the 
smallest {\sf Considered} value is always at the root of the tree. At each step of the main body of the OLIM, a {\sf Considered} 
mesh point $\mx_{{\rm new}}$
with the smallest tentative value becomes {\sf Accepted Front}. Then the hierarchical update procedure proposed in \cite{DM}
and further developed in \cite{YPC} is implemented. It consists of two substeps. First, for all {\sf Considered} 
points in $\mathcal{N}^K_{{\rm far}}(\mx_{{\rm new}})$ proposed update values involving $\mx_{{\rm new}}$ are computed. 
Second, each {\sf Unknown} point $\mx$ in  $\mathcal{N}_{{\rm near}}(\mx_{{\rm new}})$ becomes {\sf Considered} 
and a tentative value at $\mx$ is computed using the {\sf Accepted Front} points in $\mathcal{N}^K_{{\rm far}}(\mx)$ .
This algorithm is summarized in the pseudocode below. The details of each step are elaborated in \cite{YPC}. 
\noindent\makebox[\linewidth]{\rule{\textwidth}{0.4pt}}
 \begin{algorithm}[H]
 \label{alg:HU}
\KwInitialization{
Start with all mesh points being {\sf Unknown}. Set values of $U$ at them to $\infty$.
Let $\mx^{\ast}$ be an asymptotically stable equilibrium  
located at a mesh point. 
Compute tentative values of $U$ at the points $\mx\in \mathcal{N}_{{\rm near}}(\mx^{\ast})$
and change their status to {\sf Considered}. 
}

\KwTheMainBody\\
\While { the boundary of the mesh has not been reached {\bf and} the set of {\sf Considered} points is not empty}{
 {\bf 1:} Change the status of the  {\sf Considered} point $\mx_{\rm new}$ with the smallest tentative value of $U$  to {\sf Accepted Front}. \\
{\bf 2:} Change the status of all  {\sf Accepted Front} points in $\mathcal{N}_{{\rm near}}(\mx_{\rm new})$ that no longer have {\sf Considered} points 
in their  $\mathcal{N}_{{\rm near}}$-neighborhoods to {\sf Accepted}.\\
{\bf 3:} Update all  {\sf Considered} points $\mx\in\mathcal{N}_{{\rm far}}^K(\mx_{\rm new})$. The updates must involve $\mx_{\rm new}$.\\
{\bf 4:} Change the status of each  {\sf Unknown} point $\mx\in\mathcal{N_{\rm near}}(\mx_{\rm new})$ to {\sf Considered} and update them
using the {\sf Accepted Front points} in  $\mathcal{N}_{{\rm far}}^K(\mx)$. 
}
\caption{A coarse-grained pseudocode of the OLIM. 
%The hierarchical update strategy is used in steps 3 and 4. 
%A detailed description of them is found in \cite{YPC}.
}
\end{algorithm}
\noindent\makebox[\linewidth]{\rule{\textwidth}{0.4pt}}

Now we outline the hierarchical update strategy. All details of it are worked out in \cite{YPC}. 
There are three types of updates done in the following order:
$$
\text{one-point updates $\rightarrow$ triangle updates $\rightarrow$ simplex updates}.
 $$
Let $\mx$ be a {\sf Considered} point to be updated, and  $\my \in \mathcal{N}_{{\rm far}}^K(\mx)$ 
be {\sf Accepted Front}.

{\bf One-point update.}
We connect $\mx$ and $\my$ with a line segment and approximate the geometric action \eqref{GA} 
along it using the midpoint quadrature rule $\mathcal{Q}_M(\my,\mx)$.  Then the proposed value of the quasipotential at $\mx$ is
\begin{equation}
\label{1ptu}
\mathsf{Q}_1(\my,\mx) = U(\my) + \mathcal{Q}_M(\my,\mx).
\end{equation}
If $\mathsf{Q}_1(\my,\mx) $ is less than the current tentative value $U(\mx)$, 
we replace $U(\mx)$ with it. Otherwise, we leave $U(\mx)$ unchanged.
Furthermore, we compare $\mathsf{Q}_1(\my,\mx) $ with the current minimizer of the one-point update at $\mx$ and update it 
if $\mathsf{Q}_1(\my,\mx) $ is smaller.  In step 3 of Algorithm \ref{alg:HU}, the only one-point update computed is 
$\mathsf{Q}_1(\mx_{\rm new},\mx)$.
In step 4, one-point updates are computed for all {\sf Accepted Front} points $\my\in\mathcal{N}_{{\rm far}}^K(\mx)$.

{\bf Triangle update.}
Triangle updates always involve the minimizer of the one-point update $\mx_0$. The base of an admissible triangle 
is a line segment connecting $\mx_0$ and an {\sf Accepted Front} point $\mx_1$ satisfying
$\|\myp_1-\myp_0\|_1\le 2$ and $\|\myp_1-\myp_0\|_{\infty} = 1$ where $\myp_0$ and $\myp_1$ are the indices of $\mx_0$ and $\mx_1$ respectively.
%such that the $l_2$  and $l_{\infty}$ distances between the triplets of indices of $\mx_0$ and $\mx_1$ are at most 2 and 1 respectively.
The points on the line segment $[\mx_0,\mx_1]$ are parametrized by $\lambda\in[0,1]$: $\mx_{\lambda} : = \mx_0 + \lambda(\mx_1 - \mx_0)$.
The values of $U$ on $[\mx_0,\mx_1]$ are found by linear interpolation: $U(\mx_{\lambda})\equiv U_{\lambda}: = U(\mx_0) + \lambda\left(U(\mx_1)-U(\mx_0)\right)$.
Then the triangle update is done by solving the constrained minimization problem
\begin{equation}
\label{2ptu}
\mathsf{Q}_2(\mx_0,\mx_1,\mx) = \min_{\lambda\in[0,1]}\left\{ U_{\lambda} + \mathcal{Q}_M(\mx_{\lambda},\mx)\right\}
\end{equation}
and replacing the current tentative value $U(\mx)$ with  the proposed value $\mathsf{Q}_2(\mx_0,\mx_1,\mx)$ if and only if the latter is less than the former. This replacement may take place only if an interior point solution is found. Hence, we are interested in the solution to \eqref{2ptu}
only if the minimizer $\lambda^{\ast}\in(0,1)$. Therefore, we take the derivative of the function being minimized in the right-hand side of 
\eqref{2ptu}, compare its signs at the endpoints, 
and proceed with solving the nonlinear equation only if the signs are different.

{\bf Simplex update.}
One of the vertices of the triangle at the base of an admissible simplex must be the minimizer of the one-point update $\mx_0$, and one
of its sides adjacent to $\mx_0$, let's call it $[\mx_0,\mx_1]$, must be such that the constrained minimization problem \eqref{2ptu}
has given an inner point solution $\lambda^{\ast}\in(0,1)$. 
The third vertex of the base of an admissible simplex must be an {\sf Accepted Front} point $\mx_2$
such that $l_{\infty}$ distances between the indices of $\mx_0$, $\mx_1$, and $\mx_2$ are all 1, and 
at most one of the $l_1$ distances  between their indices is 2, while the other ones are 1. 
The proposed value produced by the simplex update is the solution of the constrained minimization problem
\begin{align}
\mathsf{Q}_3(\mx_0,\mx_1,\mx_2,\mx) & = \min_{\lambda\in[0,1]}\left\{ U_{\lambda} + \mathcal{Q}_M(\mx_{\lambda},\mx)\right\}, \label{3ptu}\\
{\rm where}~~\mx_{\lambda} & =  \mx_0 + \lambda_1(\mx_1 - \mx_0) + \lambda_2(\mx_2 - \mx_0),\notag \\
U_{\lambda} & =  U(\mx_0) + \lambda_1\left(U(\mx_1)-U(\mx_0)\right) + \lambda_2\left(U(\mx_2)-U(\mx_0)\right),\notag \\
\text{subject to}~~\lambda_1&\ge 0,~~\lambda_2\ge 0,~~\lambda_1 + \lambda_2\le 1. \label{3con}
\end{align}
The warm start for solving \eqref{3ptu} is the vector 
$\lambda:=[\lambda^{\ast},0]$ where $\lambda^{\ast}$ is the minimizer of \eqref{2ptu}.
As we do it for the triangle update, we wish to quickly reject the simplex update if its minimizer is certainly lying on the boundary 
of the triangle \eqref{3con}. We use the Karush-Kuhn-Tucker (KKT) optimality conditions (see \cite{nocedal}, Chapter 12)  to do so.
They boil down (see Appendix \ref{appKKT}) to  checking whether
\begin{equation}
\label{kktcheck}
\frac{\partial}{\partial\lambda_2}\left( U_{\lambda} + \mathcal{Q}_M(\mx_{\lambda},\mx) \right) \ge 0.
\end{equation}
If \eqref{kktcheck} holds, then $[\lambda^{\ast},0]$ is a local solution to \eqref{3ptu}, and hence we reject the simplex update.
Otherwise we proceed with numerical minimization using Newton's method. If an interior point solution is found,
we replace the current tentative value $U(\mx)$ with  $\mathsf{Q}_3(\mx_0,\mx_1,\mx_2,\mx)$  provided that 
$\mathsf{Q}_3(\mx_0,\mx_1,\mx_2,\mx)<U(\mx)$. Otherwise, $U(\mx)$ remains unchanged.

We remark that the computation of the quasipotential terminates as soon as a boundary mesh point becomes {\sf Accepted Front}.
This is important because the MAP that leaves the computational domain via this point might return to it, and it is crucial 
for an accurate computation of the quasipotential that the computation follows  the MAPs.

%%%%%%
\subsection{Challenges of computing the quasipotential for stochastic Lorenz'63}
\label{sec:challenge}
\begin{figure}[htbp]
\begin{center}
\includegraphics[width=0.6\textwidth]{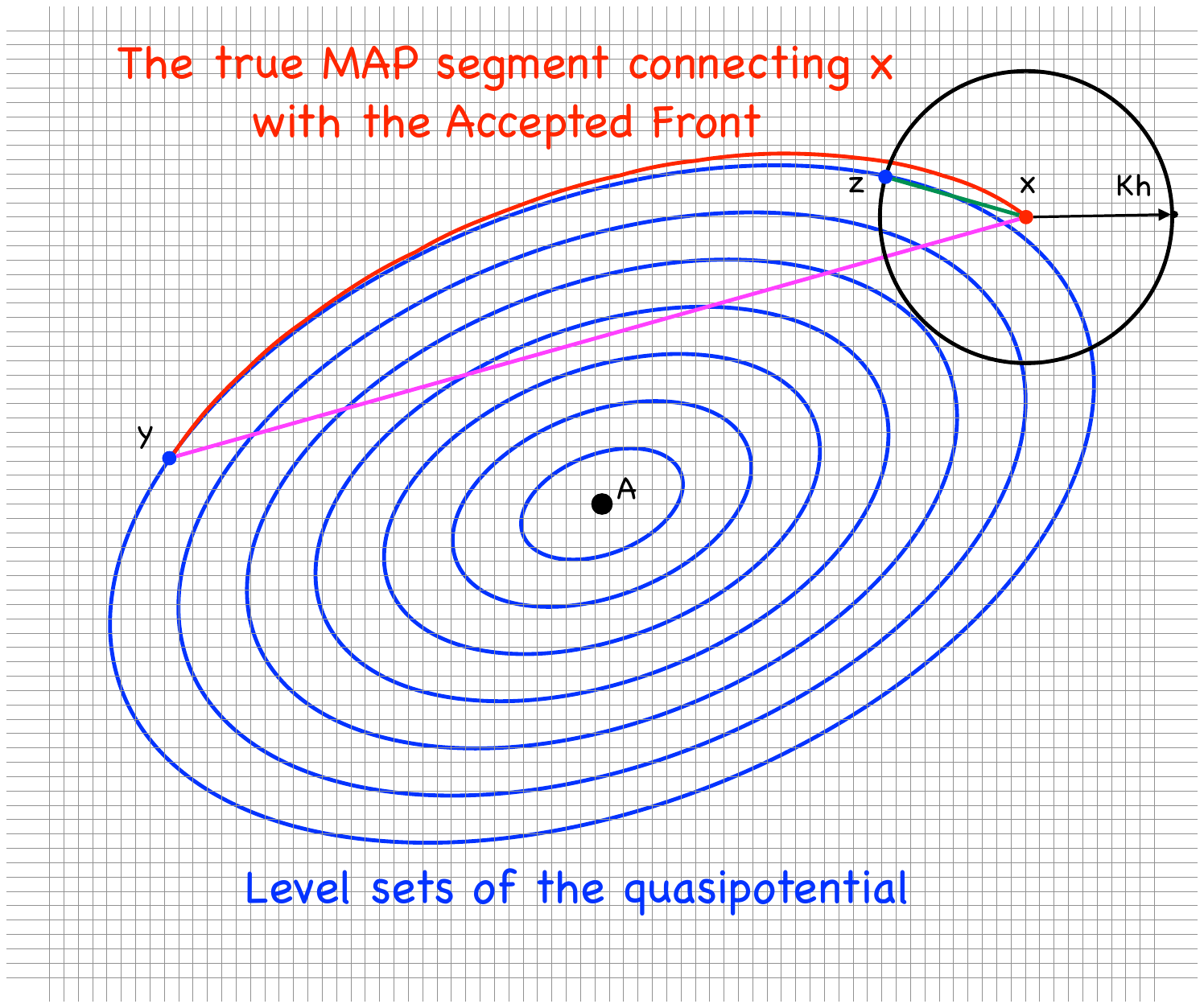}
\caption{An illustration for the difficulty of computing the quasipotential in the case where the ratio $\Xi(\mx)$ given by \eqref{Xi} is large.
The blue closed curves represent some level sets of the quasipotential. $\mx$ is a {\sf Considered} point that is up for an update.
The green segment $[\mz,\mx]$ is the best linear approximation to the MAP connecting $\mx$ with the {\sf Accepted Front} 
within the given update radius $Kh$.
}
\label{fig:challenge}
\end{center}
\end{figure}
An important characteristic of the vector field in SDE \eqref{sde0} in a neighborhood of an attractor $A$ is the 
ratio of the magnitude of the rotational component to that of the potential one \cite{YPC}:
\begin{equation}
\label{Xi}
\Xi(\mx):=\frac{\| \ml(\mx)\|}{\|\sfrac{1}{2}\nabla U(\mx)\|}.
\end{equation}
%where $U(\mx)$ is the quasipotential with respect to $A$ and $\ml(\mx): = \mb(\mx) + \sfrac{1}{2}\nabla U(\mx)$.
If $\Xi(\mx)$ is not too large (does not exceed 10) in the basin of $A$, except, perhaps some small
neighborhoods of the attractor or the escape state, the OLIMs 
 give accurate results on uniform rectangular meshes of reasonable sizes \cite{DM,DM1,YPC}.
However, if $\Xi(\mx)$ is  large (much larger that 10) in a significant part of the basin of $A$, 
the accuracy of the numerical solution by the OLIM on a regular rectangular mesh deteriorates (see Section 4 in \cite{DM}).
The problem is illustrated in Fig. \ref{fig:challenge}. Suppose the computation has reached the level set of the quasipotential
depicted with the largest closed blue curve. All mesh points inside it are either {\sf Accepted} if they have no {\sf Unknown} or {\sf Considered} nearest neighbors,
or {\sf Accepted Front}, if they do. 
Let $\mx$ be a {\sf Considered} point up for an update.
If $\Xi(\mx)$ is large, the segment of the MAP arriving at $\mx$ from the span of {\sf Accepted Front} mesh points is long.
A rough estimate for its length is $\Xi(\mx)h$ where $h$ is the mesh step. 
Let $\my$ be the point where this MAP segment starts at the span of the {\sf Accepted Front}. 
Even if the update factor $K$ were chosen large enough so that $\my$ lies in the ball centered at $\mx$ of radius $Kh$, 
the straight line segment (the magenta line segment from $\mx$ to $\my$ in Fig. \ref{fig:challenge}) 
and the midpoint quadrature rule would give poor approximations for the MAP segment and the geometric action along it respectively
resulting in an inaccurate update value at $\mx$.
It is shown in \cite{DM,YPC} that too large update factor may deteriorate the accuracy.
A safer but still too rough approximate solution would be obtained 
if the update radius is reasonably small, i.e., chosen according to the proposed rules of thumb in \cite{DM,YPC}.
Then the segment of MAP would be approximated with the green line segment $[\mz,\mx]$ in Fig. \ref{fig:challenge}.

Now imagine the case where $\Xi(\mx)\sim 10^3$ as it is for stochastic Lorenz'63 with 
$\rho_1 < \rho < \rho_2$ where the stable equilibria and the strange attractor coexist.
3D computations on regular rectangular meshes will give a qualitative idea about the geometry of the level sets of the quasipotential, but
the found quasipotential barriers will be completely off. 

\begin{figure}[htbp]
\begin{center}
\includegraphics[width=0.96\textwidth]{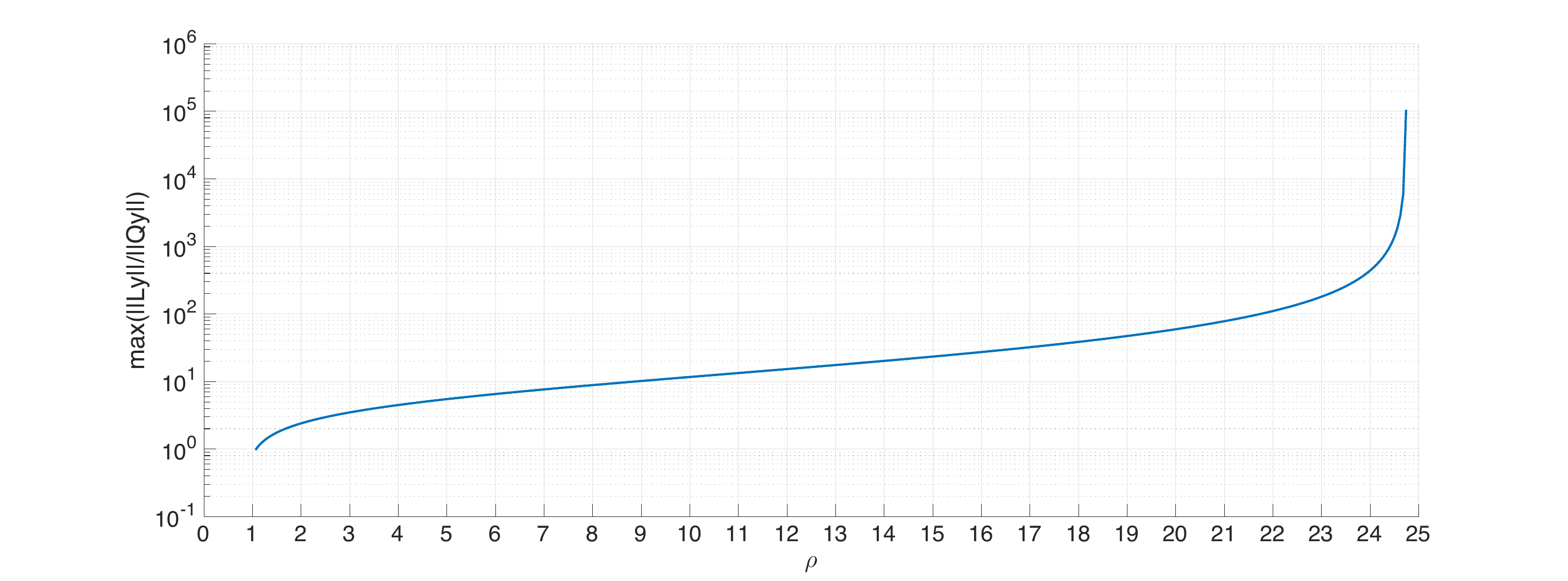}
\caption{The graph of the maximal ratio $\Xi$ of the magnitudes of the rotational and potential components of the linear SDE 
$d\my = J\my dt+\sqrt{\epsilon}d\mw$
where  $J$ is the Jacobian matrix of the right-hand side of the Lorenz system \eqref{lorenz} evaluated at the equilibrium $C_{+}$
for the range $1 < \rho< \rho_2\approx 24.74$ where $C_{+}$ is asymptotically stable.
}
\label{fig:rotpot}
\end{center}
\end{figure}

The ratio $\Xi(\mx)$ for the
Lorenz system at $1<\rho <\rho_2\approx 24.74$
can be estimated from that for the linearized system at $C_{+}$ (see Appendix \ref{sec:app0}).
The graph of $\Xi$ for the linearized system is displayed in Fig. \ref{fig:rotpot}.
It shows that the maximum of $\Xi(\mx)$ blows up as $\rho\rightarrow\rho_2$.
At $\rho = 24.4$, the largest $\rho$ at which we present the results of our computations, the maximal value of $\Xi(\mx)$ for the linearized system is 973.4.

Challenged by this problem, we have developed an approach that allows us to obtain reasonably accurate  values of the quasipotential barriers.
It consists in finding approximate 2D manifolds (or unions of 2D manifolds) where the MAPs emanating from the attractor
are located, building so-called radial meshes on them, and adjusting the OLIM for performing computations on radial meshes. 
This approach is suitable for any 3D SDE where the level sets of the quasipotential are thin, 
i.e., close to some 2D manifolds (see Assumption \ref{assumption1} below), which can be determined by visual inspection
of the computed 3D level sets. 
Note that this is a safe diagnosis as the 3D OLIM tends to make the level sets thicker than the true ones if $\Xi(\mx)$ is large.
In this case, the MAP going from the attractor to the escape state 
will be very close to any 2D manifold (or union of manifolds) approximating the level set containing the escape state. 
We find such a manifold using the characteristics of the corresponding ODE. The following lemma is instrumental for this approximation.
\begin{lemma}
\label{lemma1}
Let $A$ be an attractor of $\dot{\mx} = \mb(\mx)$, where $\mb\in C^1(\mr^3)$.
Let 
$$
\mathcal{V}_a:=\{\mx\in\mr^3~|~U(\mx)\le a\}
$$ 
be a sublevel set of the quasipotential completely lying in the basin of $A$, and 
$\gamma$ be a curve lying on the boundary of $\mathcal{V}_a$, i.e., for any $\mx\in \gamma$, $U(\mx) = a$.
Let $\mathcal{M}'$ and $\mathcal{M}$ be the manifolds consisting, respectively, of the MAPs going from $A$ to $\gamma$,  and the characteristics
starting at $\gamma$ and running to $A$.
Then $\mathcal{M}'\subset\mathcal{V}_a$ and $\mathcal{M}\subset\mathcal{V}_a$.
%, i.e., the manifolds $\mathcal{M}$ and $\mathcal{M}'$
%completely lie in $\mathcal{V}_a$.
\end{lemma}
A proof of Lemma \ref{lemma1} in found in Appendix \ref{appL1}.

Let $\gamma$ be an unstable limit cycle serving as the escape state from the basin of an attractor $A$. Let the
quasipotential at $\gamma$ be $U_{\gamma}$. We can consider a sublevel set $\mathcal{V}_a$ for
$a<U_{\gamma}$ and arbitrarily close to $U_{\gamma}$. By Lipschitz continuity of the quasipotential \cite{quasi}, 
$a$ can be chosen so that the distance between $\gamma$ and $\mathcal{V}_a$ is smaller than any given positive number.
Correspondingly, we can pick a curve $\gamma'$ lying on the boundary of $\mathcal{V}_a$ located arbitrarily close to 
the limit cycle $\gamma$. By Lemma \ref{lemma1}, 
the manifolds $\mathcal{M}'$ and $\mathcal{M}$ consisting of 
MAPs/characteristics running to/from $\gamma'$ will lie in $\mathcal{V}_a$.

\begin{assumption}
\label{assumption1}
Suppose that the level set $\mathcal{V}_a$ is close to both manifolds $\mathcal{M}$ and $\mathcal{M}'$, i.e., 
the Housdorff distances\footnotemark[3] between $\mathcal{V}_a$ and $\mathcal{M}$ and between $\mathcal{V}_a$
and  $\mathcal{M}'$ are less than some small $\delta>0$:
$$
d_H(\mathcal{V}_a,\mathcal{M})< \delta \quad {\rm and}\quad d_H(\mathcal{V}_a,\mathcal{M}')< \delta.
$$
\footnotetext[3]{$d_H(\mathcal{X},\mathcal{Y}) = \max\left\{\sup_{\mx\in\mathcal{X}}\inf_{\my\in\mathcal{Y}} \|\mx - \my\|,
\sup_{\my\in\mathcal{Y}}\inf_{\mx\in\mathcal{X}} \|\mx - \my\|\right\}$.}
\end{assumption}
Under Assumption \ref{assumption1},  the triangle inequality implies that 
the Housdorff distance between $\mathcal{M}$ and $\mathcal{M}'$ is bounded by $2\delta$:
\begin{equation}
\label{Astar}
d_H(\mathcal{M},\mathcal{M}') \le d_H(\mathcal{M},\mathcal{V}_a) + d_H(\mathcal{M}',\mathcal{V}_a)
<  2\delta.
\end{equation}

We will employ Assumption \ref{assumption1}  for $15\le \rho\le 24.4$. Figs. \ref{fig15} and \ref{fig20}
below illustrate it: compare the MAPs (the dark red curves) and the characteristics (the dark blue curves) in these figures
and observe that they lie on close manifolds located inside visibly thin level sets.

Note that the manifold $\mathcal{M}$ can be readily sampled by shooting characteristics from $\gamma'$ to $A$.
In the next section, we describe how to build radial meshes on $\mathcal{M}$, 
adjust the OLIM for them, and test its performance.

%%%%%%%%%
\subsection{Radial meshes on manifolds}
\label{sec:radial}
We call a mesh \emph{radial} if it is set up as follows.
Let $\gamma_0$ be a point or a closed curve, and let $\gamma$ be another closed curve.
We pick a finite set of simple closed curves that do not intersect pairwise and index them $\gamma_i$, $i = 1,\ldots,N_r-2$.
We add $\gamma_0$ and $\gamma_{N_r-1}\equiv \gamma$ to this set.
These curves will be referred to as \emph{parallels}.
We also pick a finite set of curves, \emph{meridians}, 
going from $\gamma_0$ to $\gamma$ and crossing each $\gamma_i$ exactly once in the order of increase of their indices.
We index the meridians from $0$ to $N_a-1$ and identify meridian 0 with meridian $N_a$. 
The resulting mesh has size $N_r\times N_a$.
Examples of radial meshes for the Lorenz system defined on manifolds consisting of all 
characteristics going from saddle cycles to asymptotically stable equilibria at $\rho = 15$ and $\rho = 24.4$ are shown in
Figs. \ref{fig2D15}(a) and \ref{fig24p4Cgamma}(a) respectively. 
A radial mesh defined between two closed curves, 
the saddle cycle $\gamma_{-}$  and a closed curve approximating an ``eye" of the strange attractor at $\rho=24.4$, is displayed in Fig. \ref{fig24p4mq}(a).
Our technique for building radial meshes is described in Appendix \ref{sec:appA} and implemented in the Matlab code {\tt make2Dmesh.m}.

To adjust the OLIM for radial meshes, we redefine the neighborhood $\mathcal{N}_{far}((i_r,i_a))$ from which a mesh
point  indexed $(i_r,i_a)$ can be updated using two update factors, radial $K_r$ and angular $K_a$, as follows:
$\mathcal{N}_{far}((i_r,i_a))$ consists of all mesh points $(j_r,j_a)$ satisfying
\begin{align*}
\max\{0,i_r - K_r\}&\le j_r\le\min\{i_r+K_r,N_r-1\}~~{\rm and }\\
|(j_a - i_a)\mod N_a|&\le K_a.
\end{align*}

Let us check whether the OLIM applied to a system with large ratio $\Xi$ produces small enough
errors on 2D radial meshes of reasonable sizes 
and these errors properly decay with mesh refinement. 
We set up an ad hoc 2D example with an asymptotically stable spiral point at the origin
and an unstable limit cycle $\|\mx \| = 1$:
\begin{equation}
\label{adhoc}
\left[\begin{array}{c}dx_1\\dx_2\end{array}\right] = \left[\begin{array}{cc} \|\mx\|^2 - 1&a \\ 
-a&  \|\mx\|^2 - 1 \end{array}\right]\left[\begin{array}{c}x_1\\x_2\end{array}\right] dt 
+ \sqrt{\epsilon}d\mw.
\end{equation}
We pick $a = 10^3$, then $\Xi \ge 10^3$.
The exact quasipotential for \eqref{adhoc} with respect to the origin is given by
\begin{equation}
\label{aexact}
U(\mx) = \begin{cases}\|\mx\|^2\left(1 -0.5\|\mx\|^2\right)&\|\mx\|\le 1 \\ 0.5,&\|\mx\|>1
\end{cases}.
\end{equation}

\begin{figure}[htbp]
\begin{center}
(a)\includegraphics[width = 0.45\textwidth]{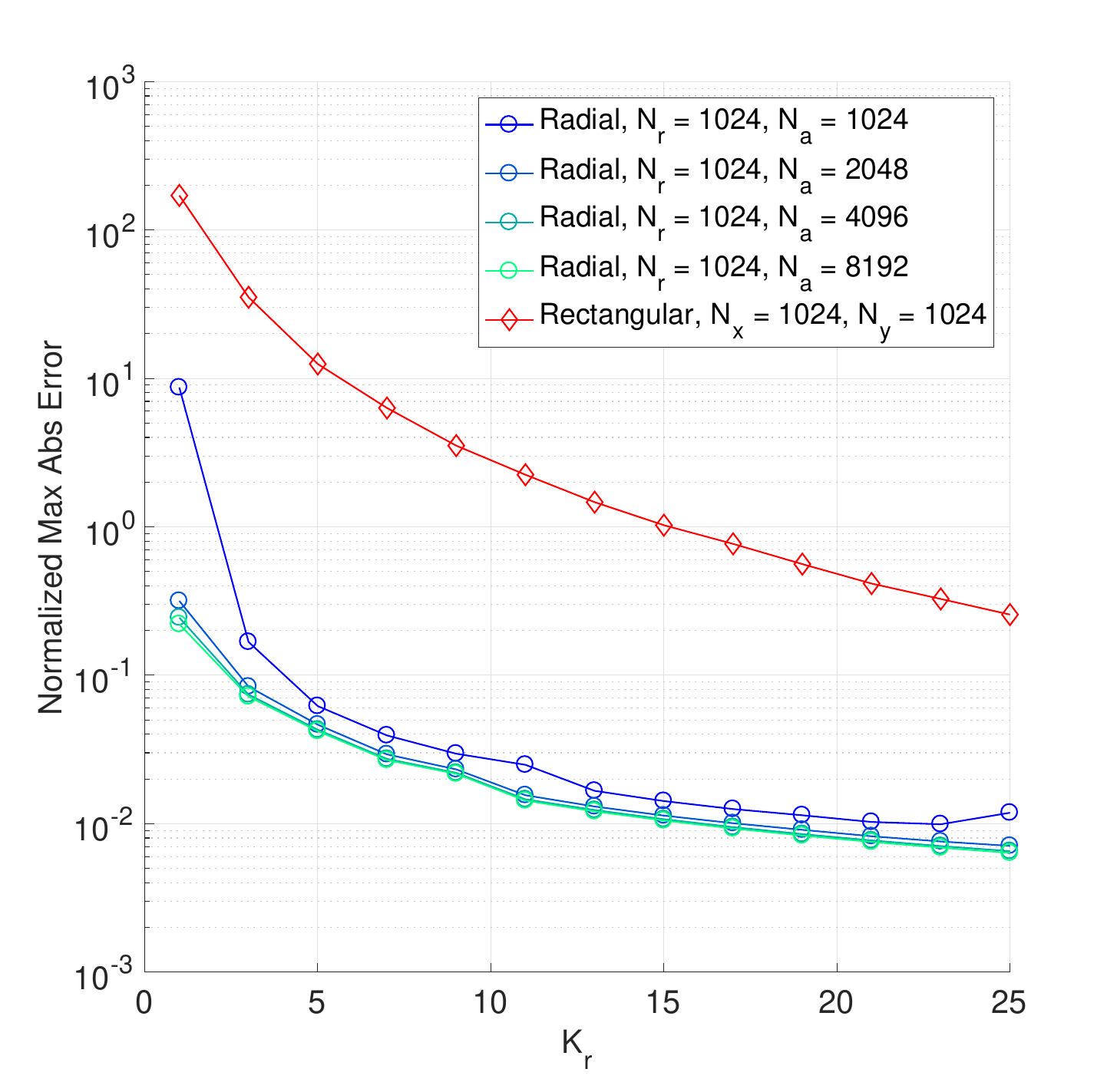}
(b)\includegraphics[width = 0.45\textwidth]{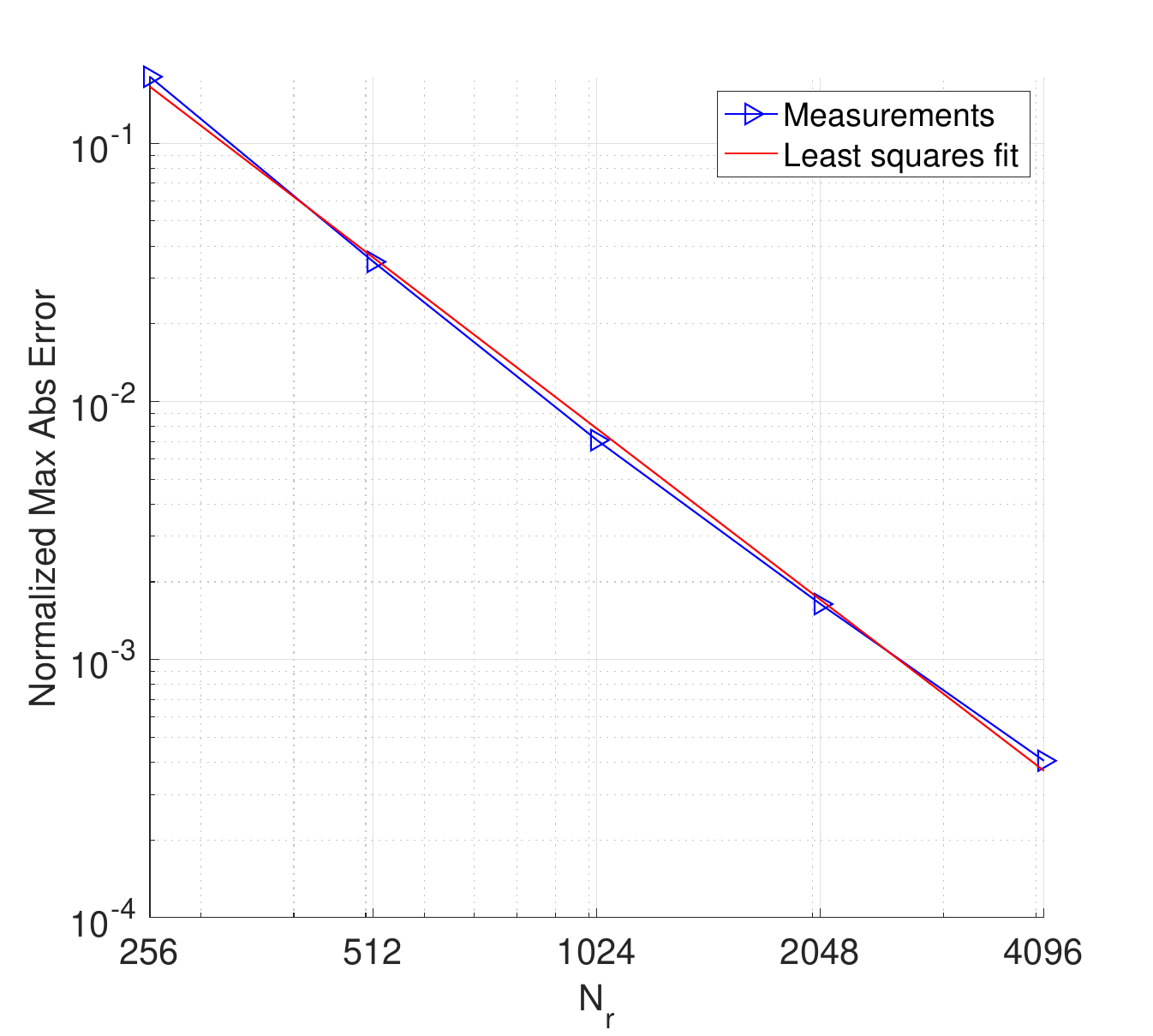}
\caption{
Measurements of numerical errors for radial meshes
$N_r\times N_a$ in computing the quasipotential for SDE \eqref{adhoc}.
(a): The dependence of the normalized maximal absolute error \eqref{nmae} on 
the update parameter $K_r$. The parameter $K_a$ was chosen so that
$N_a/N_r = K_a/K_r$. 
(b): The dependence of the normalized maximal absolute error \eqref{nmae}  
(the blue plot) on $N_r$ with $N_a = 2N_r$, $K_r = {\tt round}(N_r/40)$, and $K_a = 2K_r$.
The least squares fit \eqref{lsfit} is included for comparison. 
 }
\label{fig:adhoc}
\end{center}
\end{figure}

We have conducted two experiments with computing the quasipotential for \eqref{adhoc}.
The goal of the first experiment is to establish the dependence
of the numerical error on the relationship between $N_r$, $N_a$, $K_r$, and $K_a$.
We set $N_r=1024$ and run the solver for $N_a = 2^qN_r$, $q = 0,1,2,3$, and $K_r$ 
varying from 1 to ${\tt round}(N_r/40) = 25$ and $K_a = 2^qK_r$, respectively.
The computational domain is the unit circle. 
The dependence of the normalized maximal absolute error
\begin{equation}
\label{nmae}
E: = \frac{\max_{i_r,i_a}|U(i_r,i_a) - U_{{\rm exact}}(i_r,i_a)|}{\max_{i_r,i_a}U_{\rm exact}(i_r,i_a)}
\end{equation}
 on $K_r$  is shown in Fig. \ref{fig:adhoc}(a).
The normalized maximal absolute error (the red curve) for the $1024\times1024$ rectangular mesh defined on the square $[-1,1]^2$
 is also provided for comparison.
These results eloquently demonstrate the superiority of the 
radial meshes for computing the quasipotential in the case where the ratio $\Xi$ is large. 
Also, the choice $K_r = {\tt round}(N_r/40)$ and ${K_a = \tt round}(N_a/40) $
is reasonable and can be used as a default setting for radial meshes.

The goal of the second experiment is to verify error decay with mesh refinement.
We have run  computations with $N_r = 2^p$, $p = 8,9,10,11,12$,  $N_a = 2N_r$, $K_r = {\tt round}(N_r/40)$, and $K_a = 2K_r$. 
The plot of the normalized maximal absolute error in Fig. \ref{fig:adhoc}(b) shows the desired convergence.
The least squares fit gives a superquadratic convergence:
\begin{equation}
\label{lsfit}
E  = 3.3\cdot10^4\cdot N_r^{-2.2}.
\end{equation}

\begin{figure}[htbp]
\begin{center}
\includegraphics[width = 0.9\textwidth]{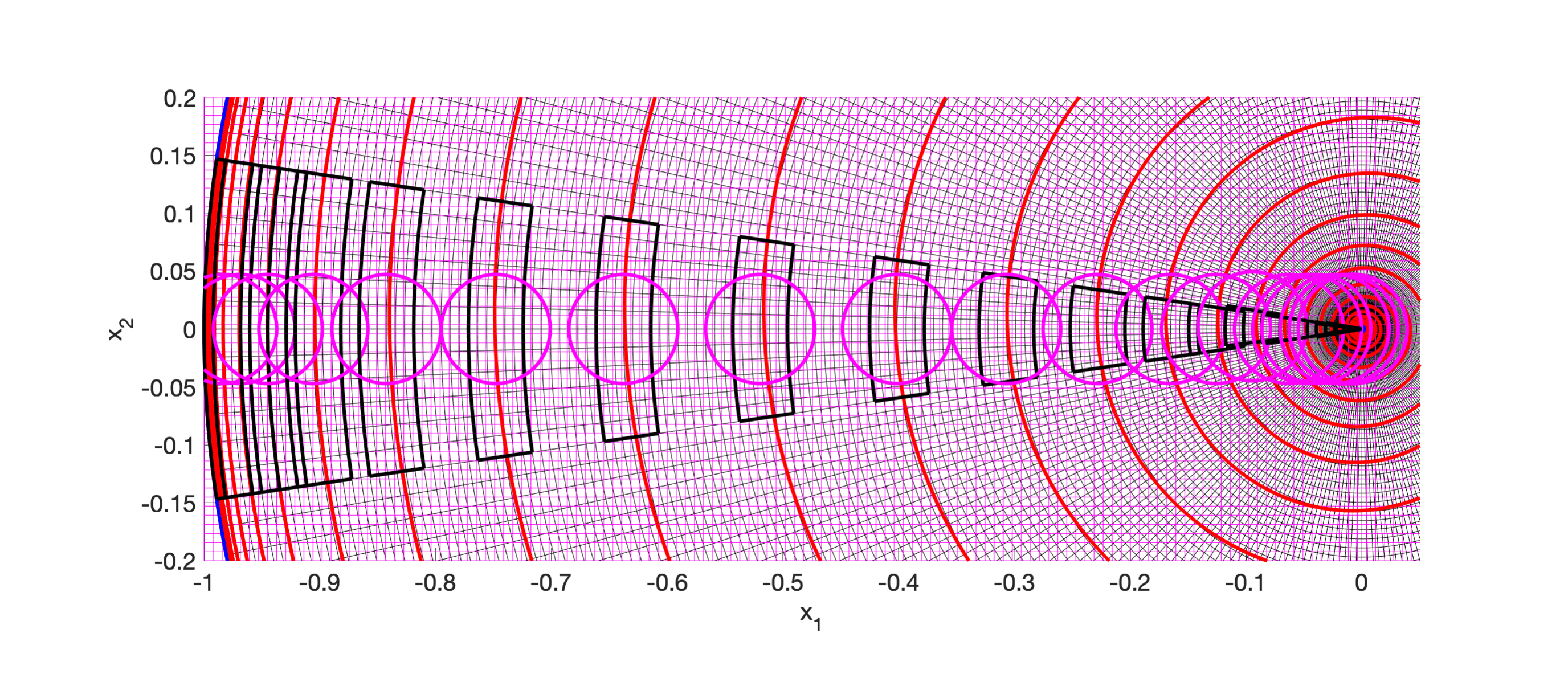}
%(b)\includegraphics[width = 0.9\textwidth]{polar_vs_rec.pdf}
\caption{An illustration explaining the advantage of radial meshes over rectangular ones for the computation of the quasipotential 
on the example
of SDE \eqref{adhoc} with $a=40$. Two computations were performed. 
The first one was done on the radial mesh with $N_r=128$, $N_a=256$, $K_r = 3$, $K_a=6$. The maximal absolute  and  RMS errors for this computations are  $1.00\cdot10^{-2}$ and $2.44\cdot10^{-3}$ respectively.
The second computation was performed on the rectangular mesh with $N = 256$ and $K =6$ and gave 
the maximal absolute  and  RMS errors of  $1.39\cdot10^{-1}$ and $6.43\cdot10^{-2}$ 
respectively which are more than an order of magnitude larger than those for the radial mesh.
The CPU times for the radial and rectangular meshes are  approximately the same: 0.24 and 0.22 seconds respectively,
The thick red curve is the exact MAP going from the equilibrium at the origin to the unstable limit cycle $r =1$. The thin black mesh is the radial mesh. The thick black curves bound some samples of its update regions. The thin magenta mesh is the rectangular mesh, and the  thick magenta circles are samples of its  update regions.
 }
\label{fig:rad_vs_rec}
\end{center}
\end{figure}

The superiority of  radial meshes over rectangular ones for the computation of the quasipotential in the basins of spiral point
 attractors of vector fields  with large rotational components is due to the fact that the radial meshes have update regions 
 better adjusted to the geometry of the MAPs than the rectangular ones. This phenomenon is illustrated
 in Fig. \ref{fig:rad_vs_rec}. The update regions of radial meshes are small near the equilibrium 
 where the MAP has high curvature and grow away from it where the MAP's curvature decreares. 
 In contrast, the update regions of rectangular meshes remain uniform. As a result, they are too large near the equilibrium and not 
 large enough away from it.

In summary, our experiments with SDE \eqref{adhoc} with a stable spiral point, an unstable limit cycle, 
and $\Xi\ge 10^3$ have demonstrated that the computation of the quasipotential on 
radial meshes of moderate sizes gives accurate and reliable results. 
%Radial meshes are advantageous for computing the quasipotential with respect to stable spiral points in 
%2D and in 3D if dimensional reduction to 2D is justified.
%Radial meshes can be upgraded to cylindric meshes if 
%the quasipotential grows at rates of the same order of magnitude in all directions.

\begin{remark}
We emphasize that we still use line segments in the OLIM on radial meshes to approximate MAP segments.
We have explored a variant of OLIM where the minimizer for each local constraint minimization problem 
is sought on the set of curves of the following form:
$$
\{(r(t),\theta(t))~|~t\in[0,1],~r(t) = r_1 + t(r_2-r_1),~\theta(t) = \theta_1+t(\theta_2-\theta_1)\}
$$
where $(r_i,\theta_i)$, $i = 1,2$, are the polar coordinates of the endpoints of the curve.
We have found that the use of line segments as in the original OLIM  gives more accurate results, so we stick with line segments.
\end{remark}

%%%%%%%%%%%%%%%%%%%%%%%%%
%%%%%%%%%%%%%%%%%%%%%%%%%

\section{Results}
\label{sec:results}
In this section, we present a collection of plots of the level sets of the computed quasipotential in 3D for the Lorenz system at
$\rho=0.5$, 12, 15, 20, and 24.4. Where appropriate, we perform 2D computations on radial meshes on manifolds and 
refine the estimates for the quasipotential barriers between different basins or regions of the phase space. 
Our collection of MAPs computed by integrating \eqref{psi} backwards in $s$ (code {\tt ShootMAPs.c}, \cite{mariakc})
can be compared with that obtained in \cite{zhou} for somewhat different set of values of $\rho$
using the minimum action method (MAM).
Note that, while the MAM is easier to program than the OLIM and it is suitable for any phase-space dimension, 
its output is biased by the initial guess for the path and hence might converge to a local minimizer
in the path-space instead of the global one. Furthermore, MAM does not allow one to visualize the level sets of the quasipotential.
Estimates for quasipotential barriers are not provided in \cite{zhou} while we do it here.

%%%%%%%%%%%%%%%%%%%%%%%%%

\subsection{$0< \rho < 1$}
\begin{figure}[htbp]
\begin{center}
(a)\includegraphics[width = 0.7\textwidth]{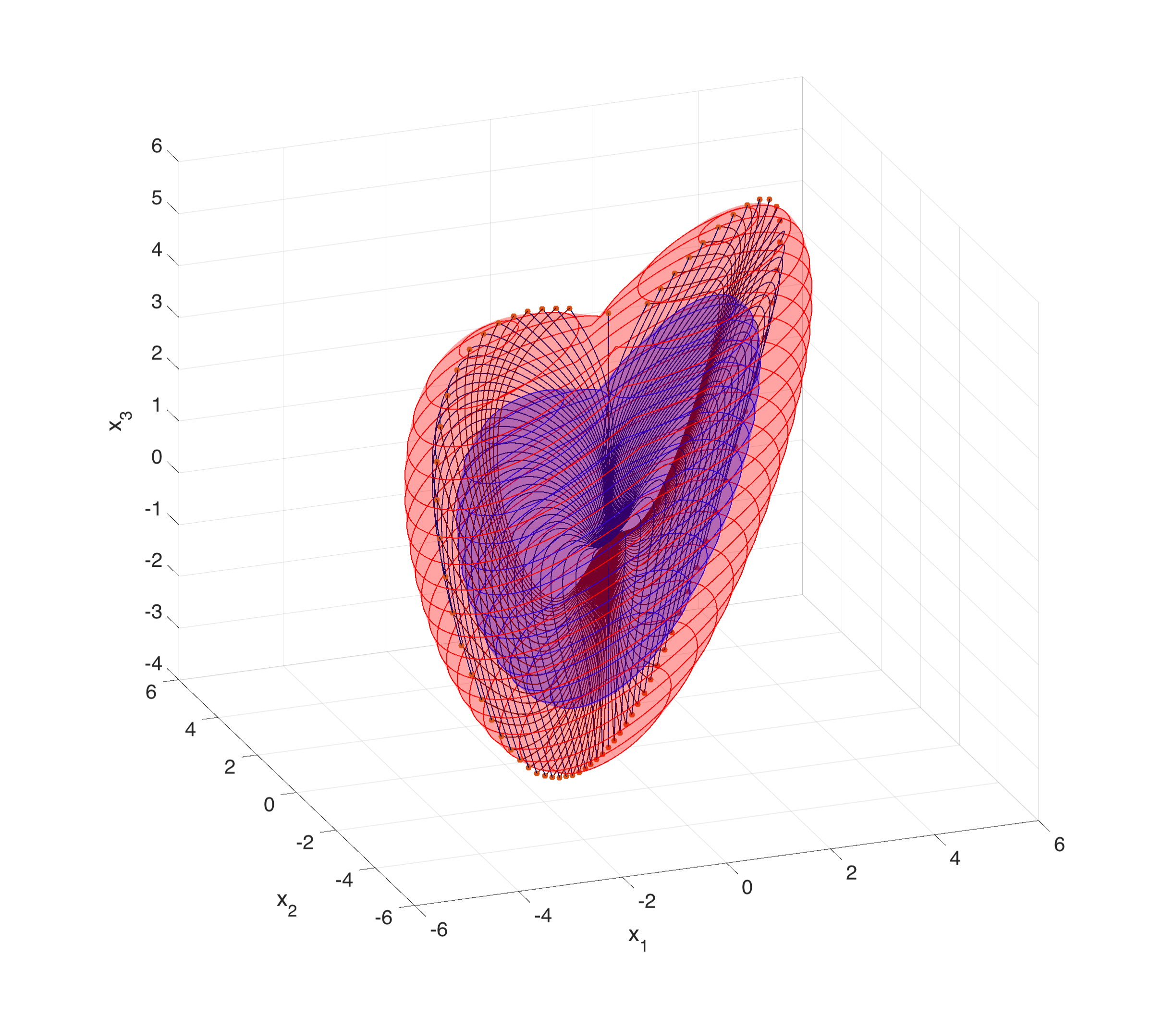}
(b)\includegraphics[width = 0.65\textwidth]{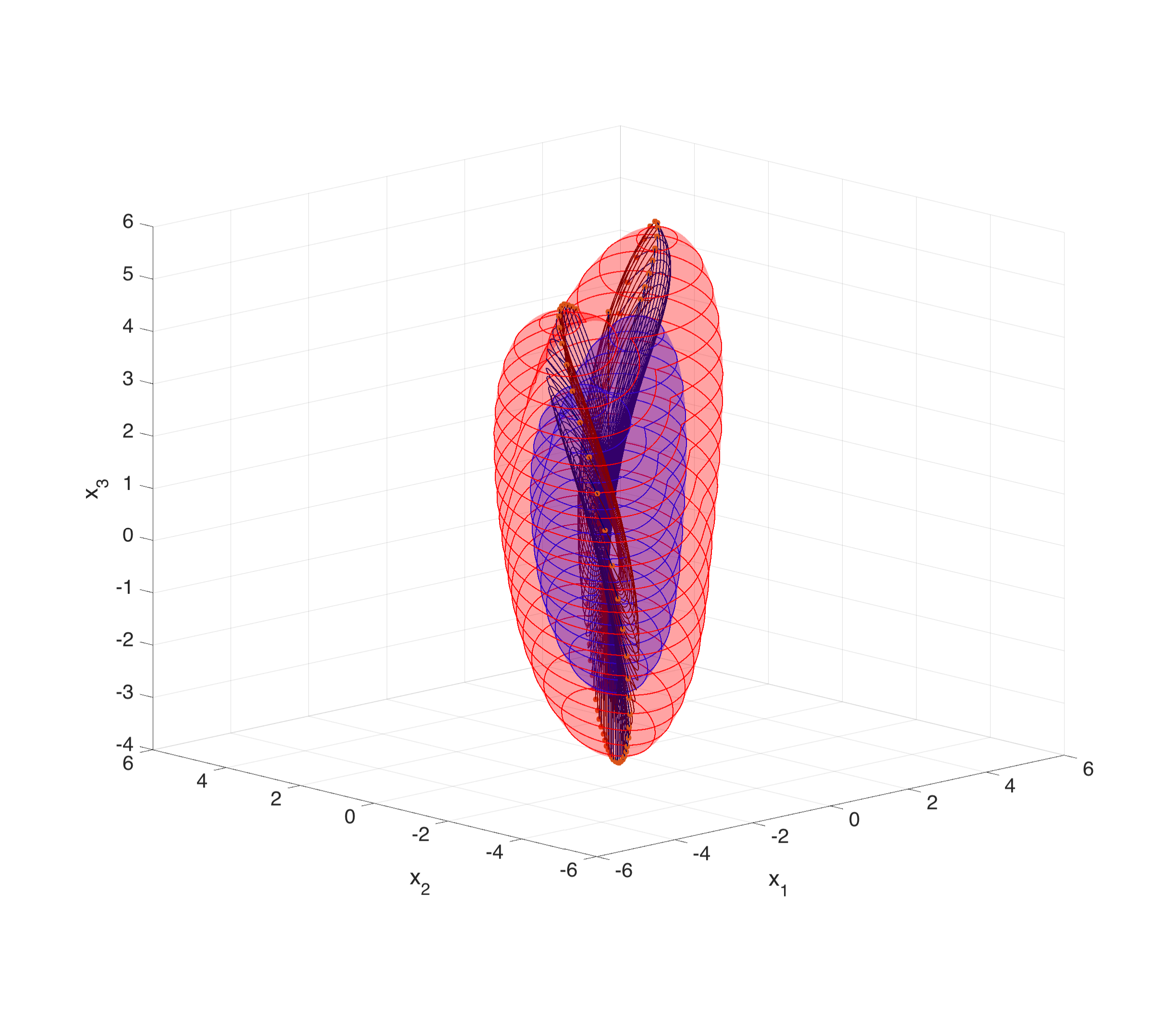}
\caption{Two views of the level sets of the quasipotential at $\rho=0.5$ corresponding to $U=20$ (the blue surface) and $U=40$ (the red surface). 
The thin blue and red closed curves lying on the corresponding level sets 
are shown to aid 3D visualization. The dark blue curves depict a collection of the
characteristics starting at the set of points marked by large orange dots
and approaching the origin.
The dark red curves represent a collection of the MAPs emanating from the origin and arriving at the same  set of points. 
A movie with this figure rotating around  the $x_3$-axis is available at \href{https://youtu.be/YscXN18lgyU}{https://youtu.be/YscXN18lgyU}.
}
\label{fig05}
\end{center}
\end{figure}

For $0< \rho < 1$, the origin is globally attracting. Two level sets of the quasipotential for $\rho = 0.5$ are shown in Fig. \ref{fig05}.
The computation was performed on $513\times513\times513$ mesh with the update factor $K = 14$. 
This choice of $K$ for $N=513$ was suggested in \cite{YPC}.
The level sets are heart-shaped and oriented approximately along the plane $x_1=x_2$.
Let $X$ be a level set and let $\gamma_X$ be the intersection of $X$ with the vertical plane $x_1=x_2$. 
The curve $\gamma_X$ runs approximately along the edge of the heart-shaped level set $X$.
We pick $X$ to be a level set corresponding to one of the largest computed values of the quasipotential and find a collection of points
marked with large orange dots lying on the corresponding curve $\gamma_X$ and forming angles from $0$ to $2\pi$ with step $\sfrac{\pi}{72}$.
The  characteristics of \eqref{lorenz} (the dark blue curves)  and the MAPs of \eqref{sde1} (the dark red curves) starting and arriving
at this set of points, respectively, are notably different. The set of characteristics starting at $\gamma_X$ and the set of MAPs
arriving at $\gamma_X$ form visibly distinct 2D manifolds.

Let us find the directions along which typical characteristics and typical MAPs approach the origin and emanate from it, respectively.
It is hard to see  in Fig. \ref{fig05} whether they coincide or not.
Let $J$ be the Jacobian matrix of the right-hand side of \eqref{lorenz} evaluated at the origin:
\begin{equation}
\label{J05}
J = \left[\begin{array}{ccc}-\sigma&\sigma&0\\ \rho&-1&0\\0&0&-\beta\end{array}\right].
\end{equation}
For the linear SDE 
\begin{equation}
\label{linsde}
d\mx = J\mx dt + \sqrt{\epsilon}d\mw,
\end{equation} 
the quasipotential decomposition is given by $J\mx = -Q\mx + L\mx$ (see Appendix \ref{sec:app0}), where
$Q$ and $L$ are matrices.
The quasipotential is the quadratic form $U(\mx) = \mx^\top Q \mx$ where $Q$  
can be found analytically \cite{quasi}:
\begin{align}
Q &= \left[\begin{array}{cc}Q_1&\\&\beta\end{array}\right],~~{\rm where} \label{Q05}\\
Q_1 &= \frac{\sigma+1}{d}\left[\begin{array}{cc} \sigma(\sigma + 1) + \rho(\rho - \sigma)&-\rho - \sigma^2 \\
-\rho - \sigma^2 & (\sigma + 1) - \sigma(\rho - \sigma)\end{array}\right],\notag\\
d &= (\sigma + 1)^2 + (\rho + \sigma)^2.
\end{align}
The rotational matrix $L = J + Q$ is
\begin{align}
L &= \left[\begin{array}{cc}L_1&\\&0\end{array}\right],~~{\rm where} \label{L05}\\
L_1 &= \frac{\rho -\sigma}{d}\left[\begin{array}{cc} \rho + \sigma^2 &- (\sigma + 1) + \sigma(\rho - \sigma)   \\
\sigma(\sigma + 1)+\rho(\rho-\sigma) & -\rho - \sigma^2  \end{array}\right].\notag
\end{align}
For the linear SDE \eqref{linsde}, MAPs are the characteristics of $\dot{\mx} = (Q + L)\mx$.
Obtaining spectral decompositions of $J = -Q + L$ and $\tilde{J} = Q + L$ for $\rho = 0.5$ we find that typical characteristics of \eqref{lorenz}
approach the origin tangent to the line ${\rm span}(\mathbf{v})$, while typical MAPs emanate from the origin 
tangent to the line ${\rm span}(\tilde{\mathbf{v}})$, where
\begin{equation}
\label{vv}
\mathbf{v} \approx \left[\begin{array}{c} 0.7241 \\  0.6897 \\ 0\end{array}\right],~~
\tilde{\mathbf{v}} \approx \left[ \begin{array}{c} 0.6924 \\ 0.7215 \\ 0\end{array}\right].
\end{equation}

%%%%%%%%%%%%%%%%%%%%%%%%%

\subsection{$1< \rho < \rho_0 \approx 13.926$}
\begin{figure}[htbp]
\begin{center}
\includegraphics[width = 0.75\textwidth]{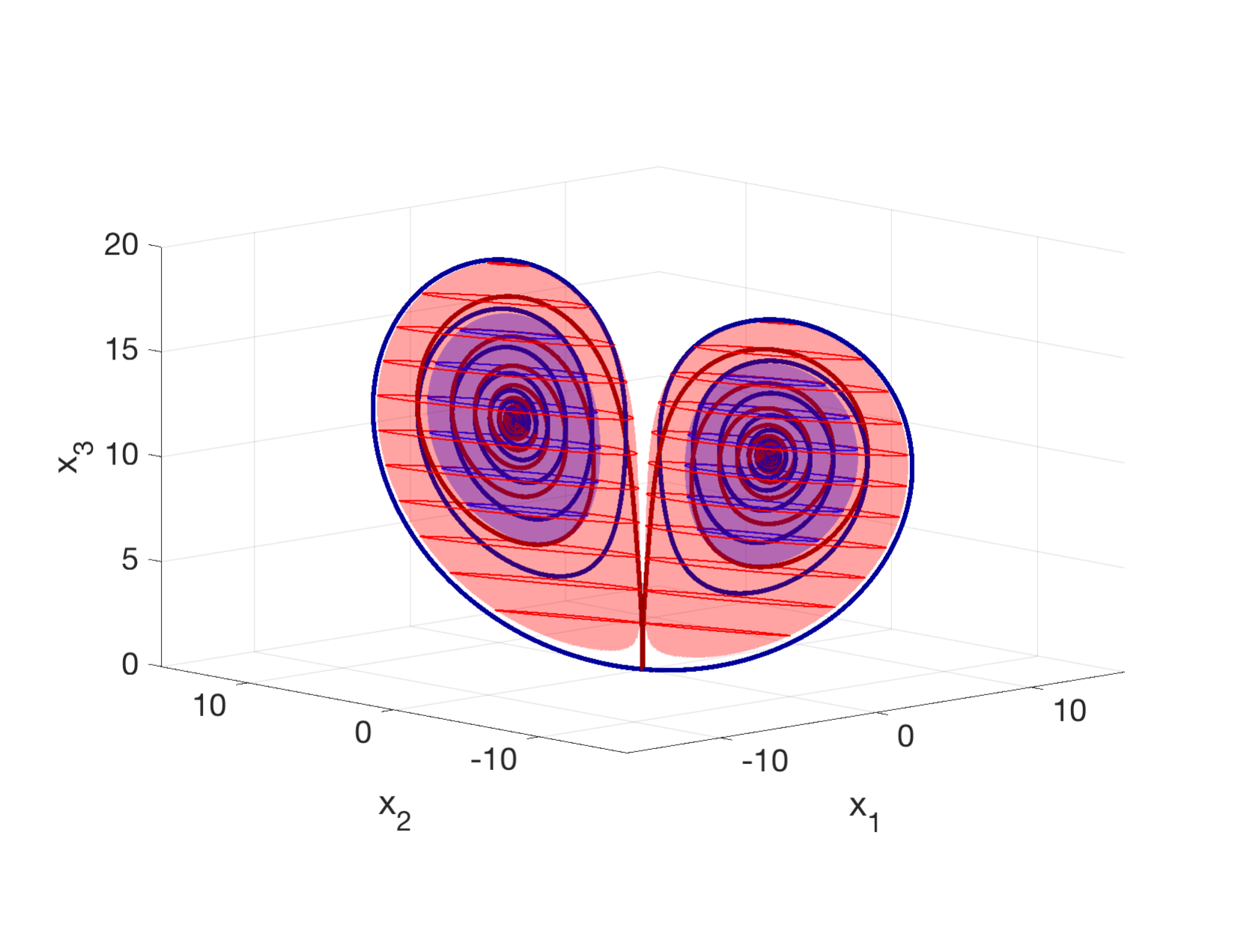}
%(b)\includegraphics[width = 0.70\textwidth]{rho12view2.pdf}
\caption{The level sets of the quasipotential at $\rho=12$ 
corresponding to $U=10$ (the blue surface) and $U=19.42$ (the red surface). 
The dark blue curves are the characteristics emanating from the origin along  its unstable directions $\pm\mathbf{\xi}$ \eqref{uvec}
and arriving at $C_{\pm}$ respectively.
The dark red curves are the MAPs going from $C_{\pm}$ to the origin.
The MAP from $C_{\pm}$ to $C_{\mp}$ is obtained by the concatenation of the MAP from $C_{\pm}$ to the origin (a dark red curve)
and the characteristic from the origin to $C_{\mp}$ (a dark blue curve).
A movie with this figure rotating around the $x_3$-axis is available at \href{https://youtu.be/-ABbuD8oDjI}{https://youtu.be/-ABbuD8oDjI}.
%\href{https://youtu.be/-ABbuD8oDjI}{here}
%%\href{https://www.dropbox.com/s/nirwjdq4vuwkaf0/LorenzQpot_rho12.00.mp4?dl=0}{here} 
%to see a movie with this figure rotating around $x_3$-axis
}
\label{fig12}
\end{center}
\end{figure}

In this interval, the equilibria $C_{\pm}$ switch from stable nodes to stable spiral points at $\rho\approx 2.1546$.
Fig. \ref{fig12} displays the level sets of the quasi-potential for $\rho = 12$ with respect to each stable equilibrium.
It was computed
on a $513\times513\times513$ mesh with $K = 14$.
The found value of the quasipotential at the origin that
serves as the transition state between $C_{\pm}$ is $19.47$.  
Therefore, at $\rho = 12$, the expected escape time from the basin of  $C_{+}$  scales as
\begin{equation}
\label{t12}
\mathbb{E}[\tau_{C_{+}}]\asymp e^{19.47/\epsilon}.
\end{equation}
The MAP from $C_{+}$ to $C_{-}$
is obtained by the concatenation of the computed MAP from $C_{+}$ to the origin (the dark red curve starting at $C_{+}$) 
and the characteristic from the origin to $C_{-}$ (the dark blue curve ending at $C_{-}$). 
Fig. \ref{fig12}(b) shows that the MAPs and the characteristics  connecting $C_{\pm}$ and the origin 
lie on close 2D manifolds.  

We did a consistency check by finding the quasipotential barrier by integrating the geometric action \eqref{GA}-\eqref{qpot1} along the found MAP and got the value
19.89 that is at a reasonable agreement with 19.47 found by our 3D computation.

%%%%%%%%%%%%%%%%%%%%%%%%%

\subsection{$13.926\approx \rho_0 < \rho < \rho_1 \approx 24.06$}
\label{sec:preturb}
In this range, the escape states from $C_{+}$ and $C_{-}$ 
are the saddle limit cycles $\gamma_{+}$ and $\gamma_{-}$ respectively.
We have computed the quasipotential for two values of $\rho$: $\rho = 15$ and $\rho = 20$.

\subsubsection{$\rho = 15$} 
\label{sec:15}
\begin{figure}[htbp]
\begin{center}
\includegraphics[width = 0.7\textwidth]{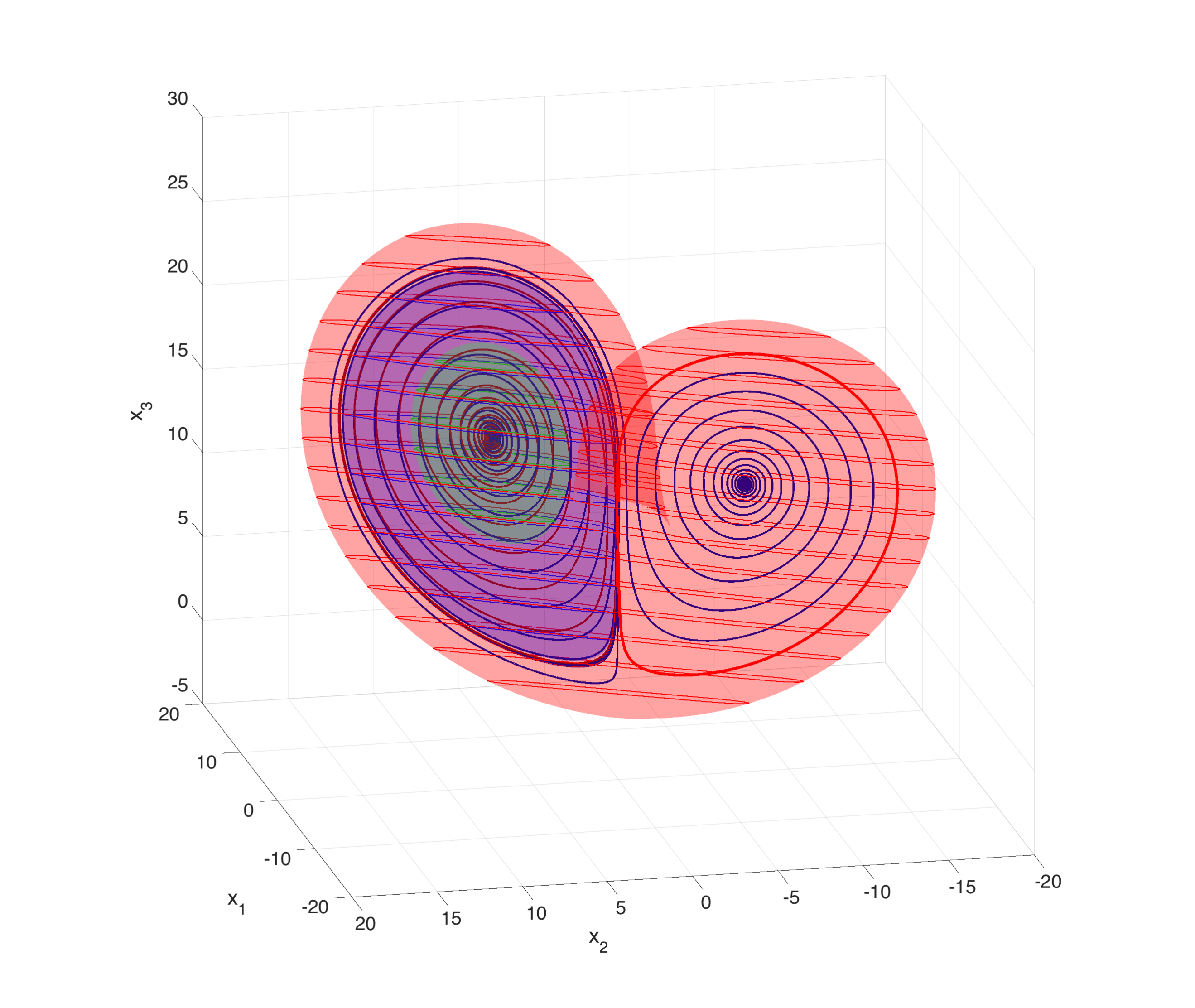}
%(b)\includegraphics[width = 0.65\textwidth]{rho15view2.pdf}
\caption{The level sets of the quasipotential at $\rho=15$ corresponding to $U=8$ (the green surface),
$U = 17.37$ (the blue surface), and $U=20$ (the red surface). 
The thick bright red curves are the saddle cycles $\gamma_{\pm}$.
The dark blue curves are characteristics running from $\gamma_{+}$
and approaching $C_{\pm}$.
The dark red curve is a MAP starting at $C_{+}$ and approaching $\gamma_{+}$. 
%Click  \href{https://youtu.be/mzdUD-ngqYs}{here} to see a movie with this figure rotating around $x_3$-axis
A movie with this figure rotating around the $x_3$-axis is available at \href{https://youtu.be/mzdUD-ngqYs}{https://youtu.be/mzdUD-ngqYs}.
}
\label{fig15}
\end{center}
\end{figure}

The computed quasipotential  for $\rho = 15$ with respect to  $C_{+}$ is visualized in Fig. \ref{fig15}.
First, we picked a large computational domain to embrace the level set of the quasipotential enclosing 
both of the stable equilibria $C_{\pm}$ and used 
a $613\times613\times613$ mesh and $K = 15$.
Second, we chose a smaller domain just to enclose $\gamma_{+}$. 
It was a cube with side 13 centered at $C_{+}$, and  the mesh in it was $1001\times1001\times1001$.
$K$ was set to 20.
The found quasipotential is nearly constant on $\gamma_{+}$: it varies between 17.42 and 17.45.
The saddle cycles $\gamma_{\pm}$ are depicted with thick bright red curves.
A maximum likelihood transition path from $C_{+}$ to $C_{-}$ can be obtained by the
concatenation of a MAP from $C_{+}$ to $\gamma_{+}$, the saddle cycle $\gamma_{+}$, and a characteristic going from $\gamma_{+}$ to $C_{-}$.
One such  MAP and one such characteristic are the dark red and dark blue curves in Fig. \ref{fig15} respectively.

\begin{figure}[htbp]
\begin{center}
(a)\includegraphics[width = 0.8\textwidth]{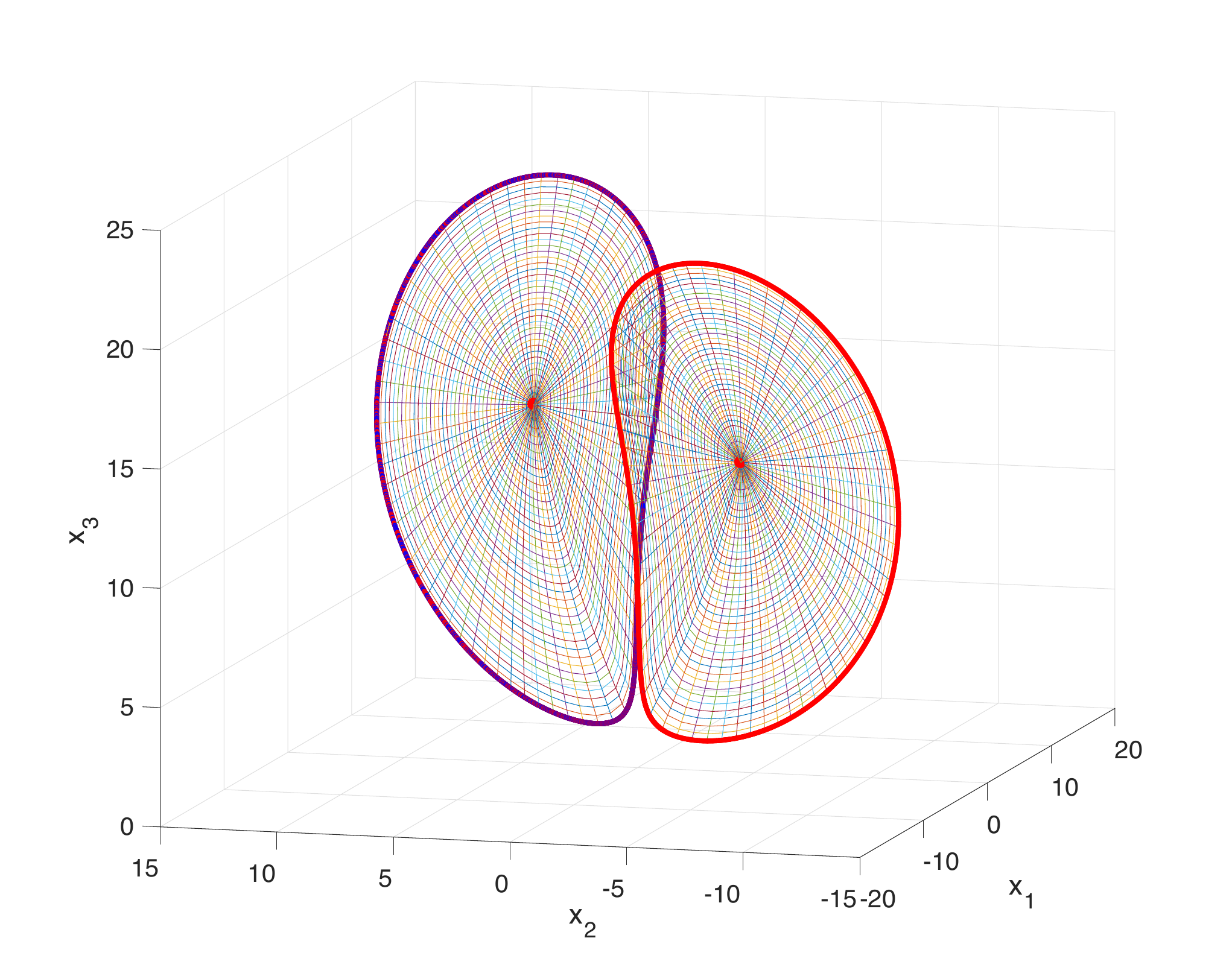}
(b)\includegraphics[width = 0.8\textwidth]{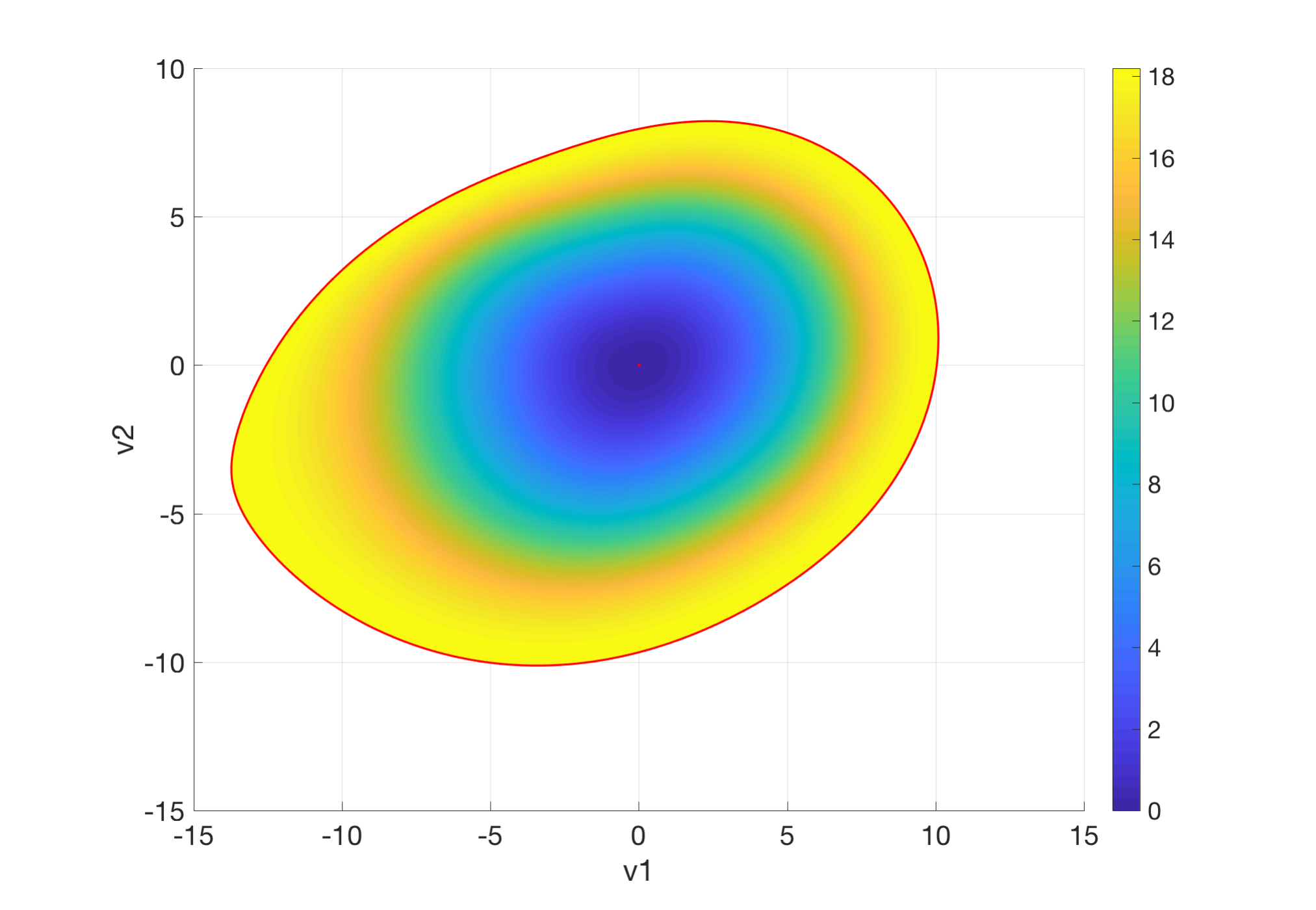}
\caption{(a): Radial meshes on the manifold $\mathcal{M}_{\pm}$
consisting of all characteristics going from the saddle cycles 
$\gamma_{+}$ (the thick purple curve) and $\gamma_{-}$ (the thick red curve) 
to the equilibria $C_{\pm}$  (the large red dots) respectively.
(b): The quasipotential computed on $\mathcal{M}_{+}$.
}
\label{fig2D15}
\end{center}
\end{figure}

Willing to refine our relatively rough 3D computation and find a more accurate value of the quasipotential on $\gamma_{+}$
with respect to $C_{+}$, we perform 2D computations on the manifold $\mathcal{M}_{+}$ consisting of all characteristics 
going from $\gamma_{+}$ to $C_{+}$ using the code \\
{\tt olim2DEquilibLimitCycle.c}. 
Fig. \ref{fig15} suggests that $\mathcal{M}_{+}$ is close to the 2D manifold
consisting of all MAPs from $C_{+}$ to $\gamma_{+}$. So, we neglect the discrepancy between them.
We generate  2D radial meshes on $\mathcal{M}_{+}$  (see Appendix \ref{sec:appA}) whose coarsened version is shown in Fig. \ref{fig2D15}(a).
The computed quasipotential on $\mathcal{M}_{+}$ is shown in Fig. \ref{fig2D15}(b).
We first ran the OLIM on a radial mesh of size $2001\times7200$ and then repeated the computation on a  refined mesh  of size $4001\times14400$. 
The radial update factors $K_r$ were $50$
and $100$ respectively, and the angular update factors $K_a$ were $180$ and $360$ respectively.
For the coarser mesh, the resulting values of the quasipotential on $\gamma_{+}$ varied from $18.19488$ to $18.19501$ averaging at $18.19495$.
For the finer mesh, these numbers were, respectively, $18.19536$, $18.19541$, and $18.19536$.
These results suggest the following estimate for expected escape time from $C_{+}$ at $\rho = 15$:
\begin{equation}
\label{t15}
 \mathbb{E}[\tau_{C_{+}}]\asymp e^{18.2/\epsilon}.
\end{equation}
For comparison and a consistency check, we have also found the quasipotential barrier by integrating the geometric action along the MAP
going from $C_{+}$ to $\gamma_{+}$. 
%The MAP is obtained by integrating \eqref{psi} backwards. 
Note that the length of this MAP is infinite. However, the contribution to the geometric action
from the integration along its infinite  piece 
lying within an $\delta$-tube around $\gamma_{+}$ tends to zero as $\delta\rightarrow 0$ as the quasipotential is Lipschitz-continuous \cite{quasi}.
Therefore, it suffices to take a finite piece of the MAP starting at $C_{+}$ and ending near $\gamma_{+}$.
We took a piece of MAP of length $308.7$ and obtained the value of the quasipotential barrier 19.3 
which is closer to 18.2 found by the 2D computation rather than to 17.4
found by the 3D one. The result 19.3 is affected by numerical errors in the MAP and by the quadrature error 
amplified by the large length of the MAP. As $\rho$ increases to $\rho_2\approx24.74$, the MAP spirals denser and denser, 
and integration of the geometric action along it becomes less and less accurate. 
So, we abandon this consistency check for values of $\rho$ larger than 15.

\subsubsection{$\rho = 20$}
\begin{figure}[htbp]
\begin{center}
\includegraphics[width = 0.8\textwidth]{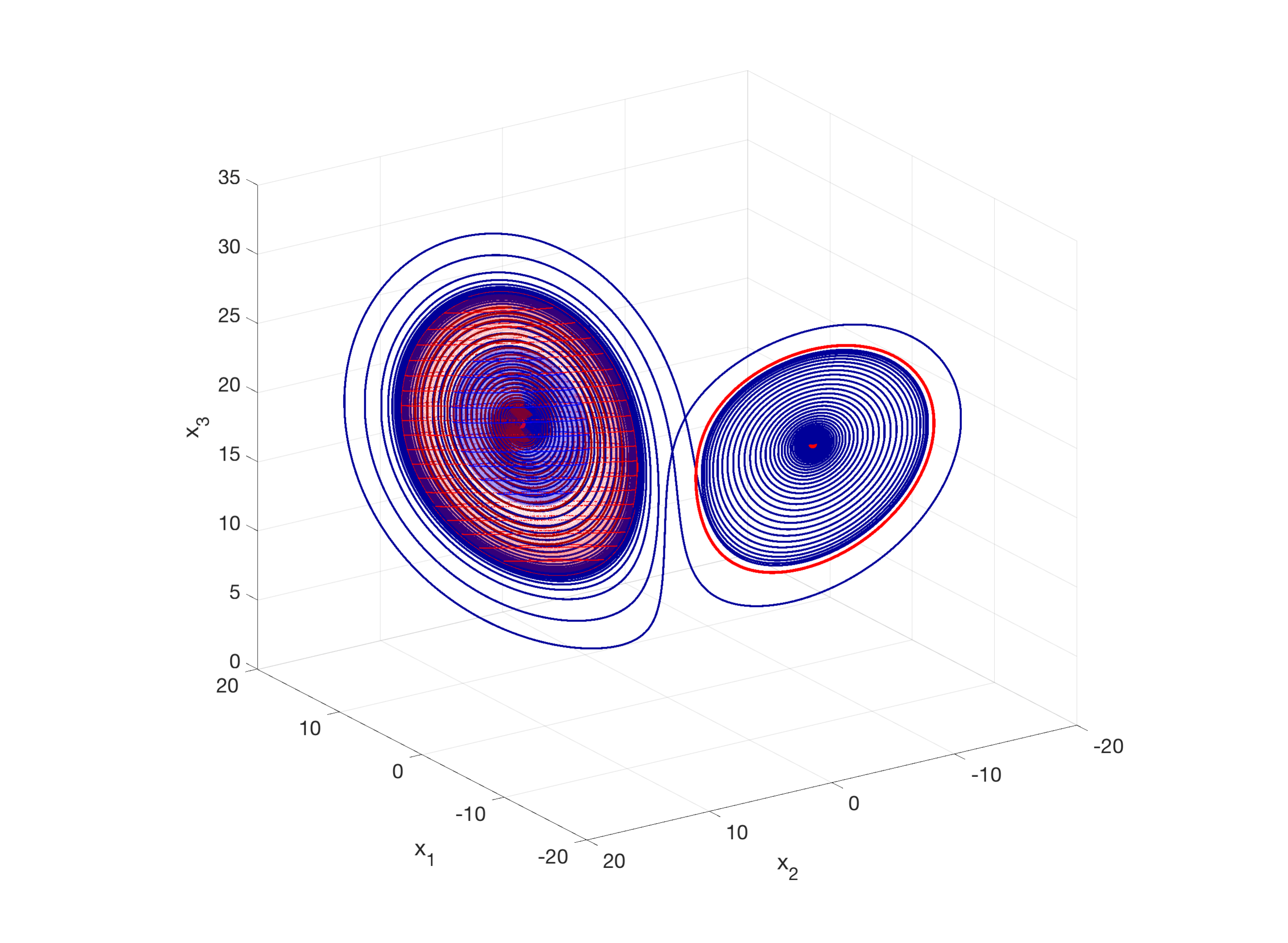}
%(b)\includegraphics[width = 0.75\textwidth]{rho20view2.pdf}
\caption{The level sets of the quasipotential at $\rho=20$ corresponding to $U=3.3$ (the blue surface),
and $U = 6.58$ (the red surface). 
The thick bright red curves are the saddle cycles $\gamma_{\pm}$.
The dark blue curves are characteristics going from $\gamma_{+}$
to $C_{+}$ and $C_{-}$.
The dark red curve is a MAP starting at $C_{+}$ and approaching $\gamma_{+}$. 
%Click  \href{https://youtu.be/JhBU0-dnos8}{here} to see a movie with this figure rotating around $x_3$-axis
A movie with this figure rotating around the $x_3$-axis is available at \href{https://youtu.be/JhBU0-dnos8}{https://youtu.be/JhBU0-dnos8}.
}
\label{fig20}
\end{center}
\end{figure}

\begin{figure}[htbp]
\begin{center}
(a)\includegraphics[width = 0.45\textwidth]{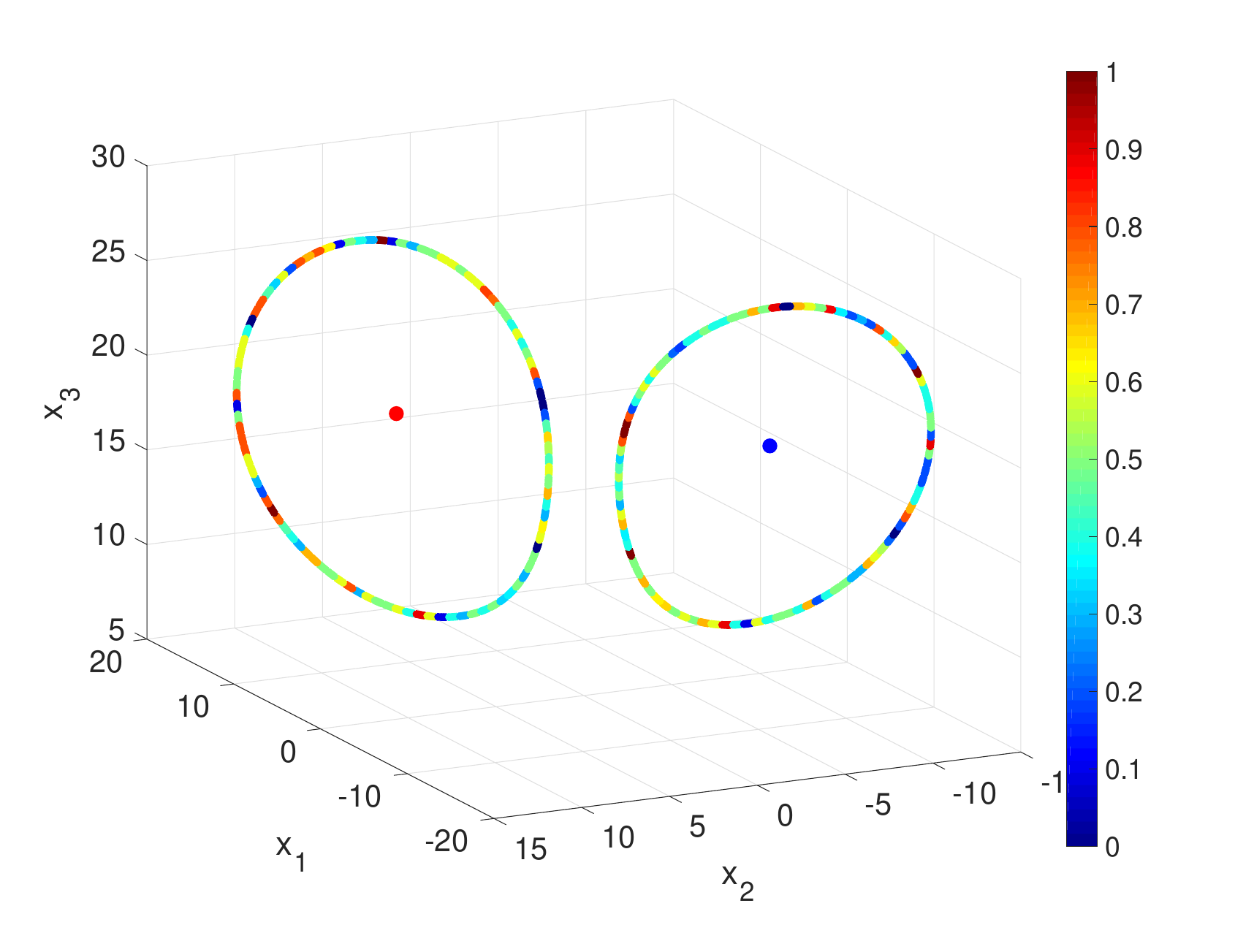}
(b)\includegraphics[width = 0.45\textwidth]{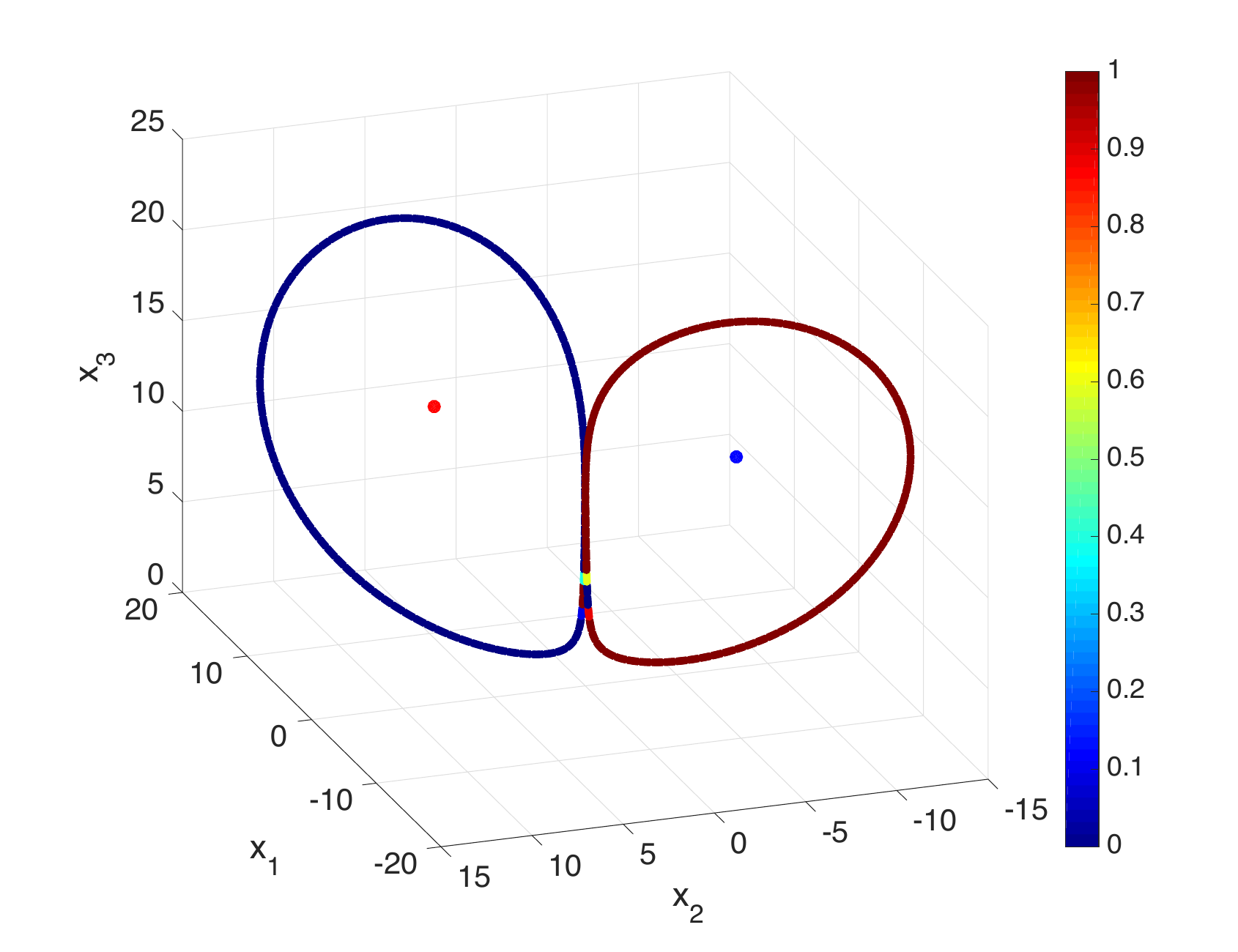}
\caption{The probability for a trajectory starting on the cones $\Upsilon_{\pm}$ at the point of the form 
$\mx + 0.002(\mx-C_{\pm})$ where $\mx\in\gamma_{\pm}$ respectively to converge to $C_{+}$.
(a): $\rho = 20$. (b): $\rho = 15$.
}
\label{fig20quest}
\end{center}
\end{figure}

For $\rho = 20$, we performed a computation in the cube with side 26 centered at $C_{+}$ on a $1001\times1001\times1001$ mesh with $K = 20$. 
This cube encloses $\gamma_{+}$. 
The values of the computed quasipotential on $\gamma_{+}$ range from $6.59$ to $6.62$ and average at $6.61$.
The level sets corresponding to $U=3.3$ and $U = 6.58$ are shown in Fig. \ref{fig20}.
A 2D computation on the manifold $\mathcal{M}_{+}$ similar to the one described in Section \ref{sec:15} 
gave $U(\gamma_{+})\in[6.1172,6.1175]$ with the average at $6.1172$.
The MAP going from $C_{+}$ to $\gamma_{+}$  as well as the characteristics going from $\gamma_{+}$ to $C_{+}$ spiral
notably denser than their counterparts at $\rho = 15$, and the level sets of the quasipotential are thinner.
The saddle cycles are the escape states from the basins of $C_{\pm}$ to a chaotic region \cite{preturbulence}
where it is hard to predict for a characteristic which attractor, $C_{+}$ of $C_{-}$, it will eventually approach.
We traced $1000$ trajectories starting on the cone $\Upsilon_{+}$ (see Table \ref{table:lorenz}) at the points of the form 
$\mathbf{y}_i: = \mx_i+0.002(\mx_i-C_{+})$ where $\mx_i \in\gamma_{+}$, $i = 1,\ldots,1000$,
are equispaced, and recorded whether they converged to $C_{+}$ or $C_{-}$ as $t\rightarrow\infty$:
508 and 492 trajectories converged to $C_{+}$ and $C_{-}$ respectively.
Then we subdivided $\gamma_{+}$ into 100 intervals of equal length
and used the recorded data to estimate the probability for a trajectory starting at each $\mathbf{y}_i$ 
corresponding to $\mx_i$ in each interval to converge to $C_{+}$.
The result is shown in Fig. \ref{fig20quest}(a). The probabilities for $\gamma_{-}$ are obtained by symmetry.
Note that a similar calculation for $\rho=15$ gave the probability distribution depicted in Fig. \ref{fig20quest}(b): 975 out of 1000 trajectories starting an the analogous 
points of the cone $\Upsilon_{+}$ eventually approached $C_{-}$, while 25 returned to $C_{+}$.  
The uncertainty for where the trajectory of \eqref{sde1}
that escapes all level sets of the quasipotential not containing the saddle cycle will eventually go, 
to $C_{+}$ or to $C_{-}$, appears where the saddle cycles $\gamma_{\pm}$
come close to each other.

Summarizing our findings for $\rho=20$, we predict that the expected escape time from $C_{\pm}$ to the chaotic region scales as
\begin{equation}
\label{t20}
 \mathbb{E}[\tau_{C_{+}}]\asymp e^{6.1/\epsilon}.
\end{equation}
%Unsurprisingly, the region of uncertainty is located around the place where $\gamma_{+}$ and $\gamma_{-}$ come close to each other.

%%%%%%%%%%%%%%%%%%%%%%%%%

\subsection{$ 24.06\approx \rho_1 < \rho <\rho_2\approx 24.74$}
\label{sec:coex}
\begin{figure}[htbp]
\begin{center}
\includegraphics[width = 0.8\textwidth]{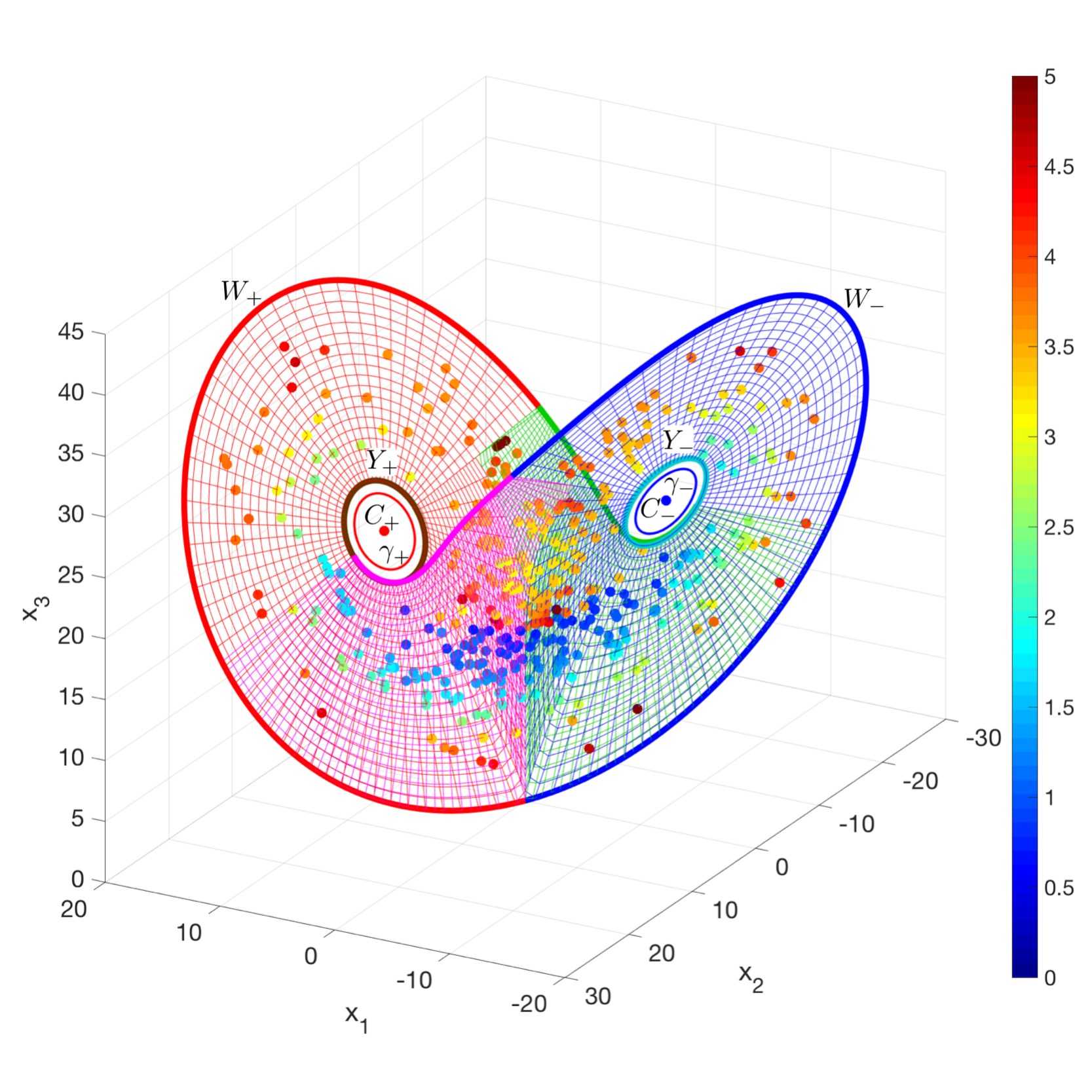}
\caption{ The strange attractor $A_L$ at $\rho = 24.4$ is approximated by a union of four manifolds: 
red, magenta, blue, and green. The color of the large dots on the manifolds indicate the thickness of the fractal (the Lorenz attractor)
at the corresponding locations. The colorbar corresponds to $-\log_{10}w(\mx)$ where $w(\mx)$ is the thickness of
the fractal near the location $\mx$. Hence dark blue dots indicate thickness $\sim10^{-1}$, light blue ones -- $\sim10^{-2}$,
yellow ones -- $\sim10^{-3}$, orange ones -- $\sim10^{-4}$, and red ones -- $\sim10^{-5}$.
}
\label{figSA}
\end{center}
\end{figure}

It was recognized by Lorenz \cite{lorenz} that the strange attractor is an ``infinite complex of surfaces'', i.e., a fractal, which is 
a very complicated geometric object. 
The addition of small white noise to the Lorenz system regularizes and simplifies its dynamics
in the sense that it renders the fine structure of the Lorenz attractor 
irrelevant and allows for a description of the dynamics in terms of probability measures. 
Taking this into account,
% for the purpose of visualization and 
%reference to different parts of it, 
we approximate the strange attractor $A_L$ with a union of four manifolds
as shown in Fig. \ref{figSA}. These manifolds were obtained using the code {\tt StrangeAttractorMesh.m} in a way
similar to the one described in Appendix \ref{sec:appA}. 
The key component of this construction is finding a trajectory going into the saddle at the origin.
 We will refer to the inner boundaries of the red and blue 
manifolds plotted with brown and cyan, respectively,  as the eyes $Y_{+}$ and $Y_{-}$.
The union of the red and green boundaries will be called wing $W_{+}$. Similarly,
the union of the blue and magenta  boundaries forms the wing $W_{-}$.
In order to understand what is the minimal reasonable value of the parameter $\epsilon$ in \eqref{sde1} that makes
such an approximation sensible, we  have estimated the thickness of the strange attractor at 398 randomly picked points.
 Details are provided in Appendix \ref{sec:appB}.
 The thickness map in Fig. \ref{figSA} indicates that the thickness of $A_L$ does not exceed $10^{-2}$ 
 wherever it is approximated by a single manifold. 
 Larger values of thickness are found in places where we approximate $A_L$ with two close manifolds. Hence they are just an
 artifact of our thickness measurement method.
 The thickness map suggests that $\sqrt{\epsilon}$ in SDE \eqref{sde1} should be at least $10^{-2}$, i.e., $\epsilon\gtrsim10^{-4}$.

\begin{figure}[htbp]
\begin{center}
\includegraphics[width = 0.7\textwidth]{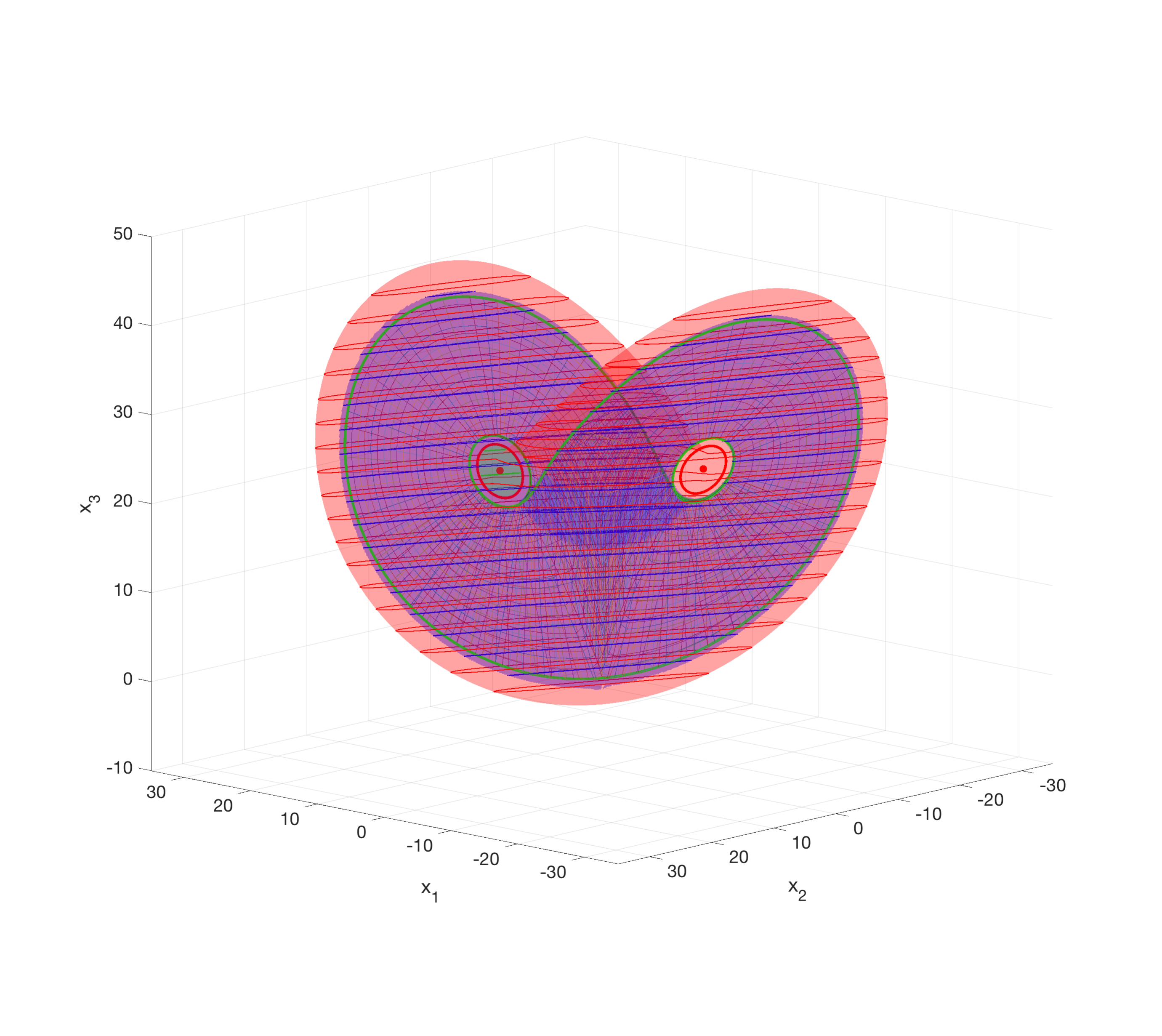}
%(b)\includegraphics[width = 0.65\textwidth]{rho24p4view2.pdf}
\caption{ $\rho = 24.4$.
The level sets of the quasipotential computed with respect to $C_{+}$. 
The green surface corresponds to the quasipotential value slightly less than the one at $\gamma_{+}$.
The blue and red ones correspond to $U = 2$ and $U = 20$ respectively.
The strange attractor is depicted with a mesh visible inside the blue and red surfaces.
A movie with this figure rotating around the $x_3$-axis is available at \href{https://youtu.be/ELqkeb8M1fg}{https://youtu.be/ELqkeb8M1fg}.
%Click  \href{https://youtu.be/ELqkeb8M1fg}{here} to see a movie with this figure rotating around $x_3$-axis
}
\label{fig24p4}
\end{center}
\end{figure}

\begin{figure}[htbp]
\begin{center}
(a)\includegraphics[width = 0.45\textwidth]{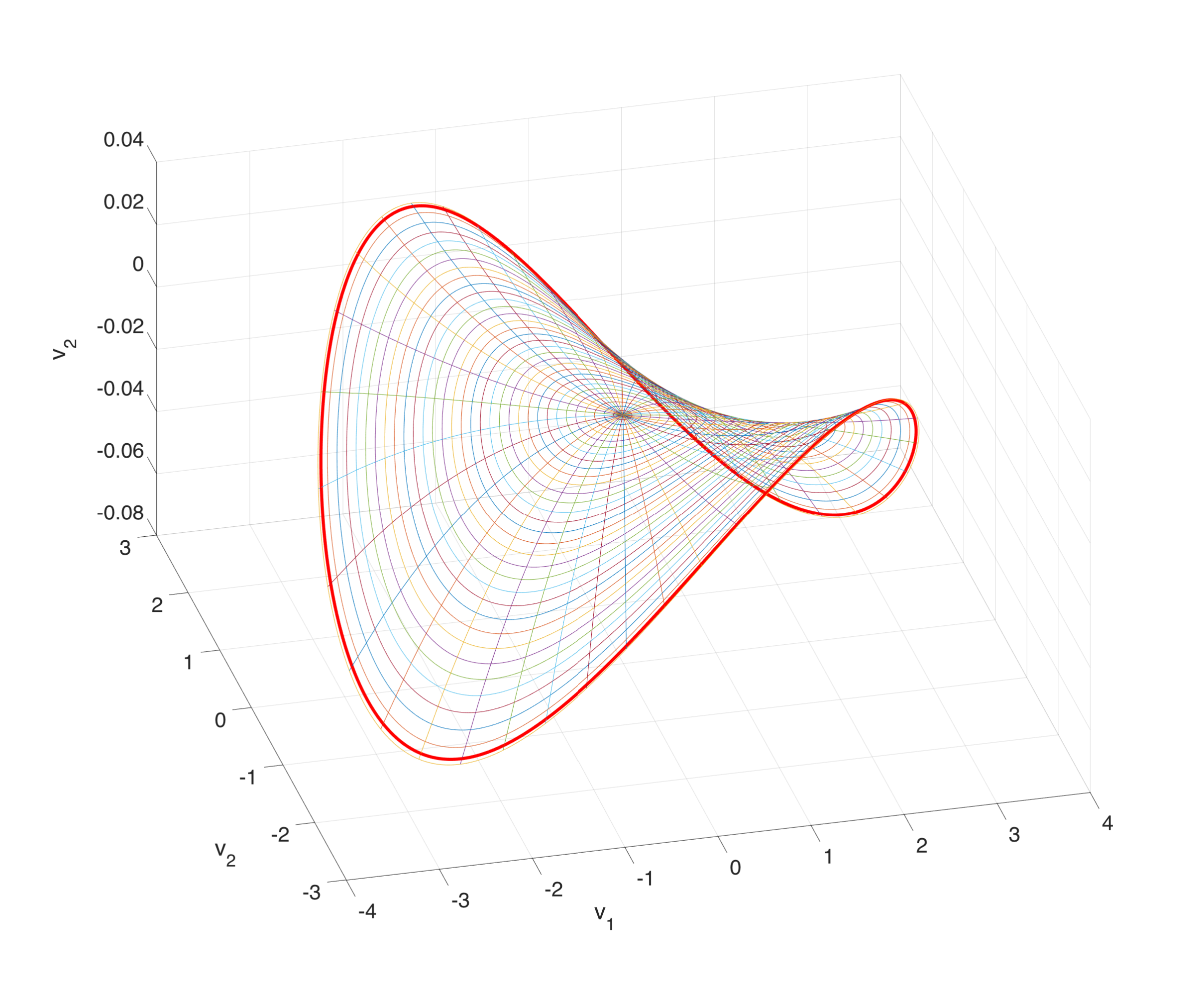}
(b)\includegraphics[width = 0.45\textwidth]{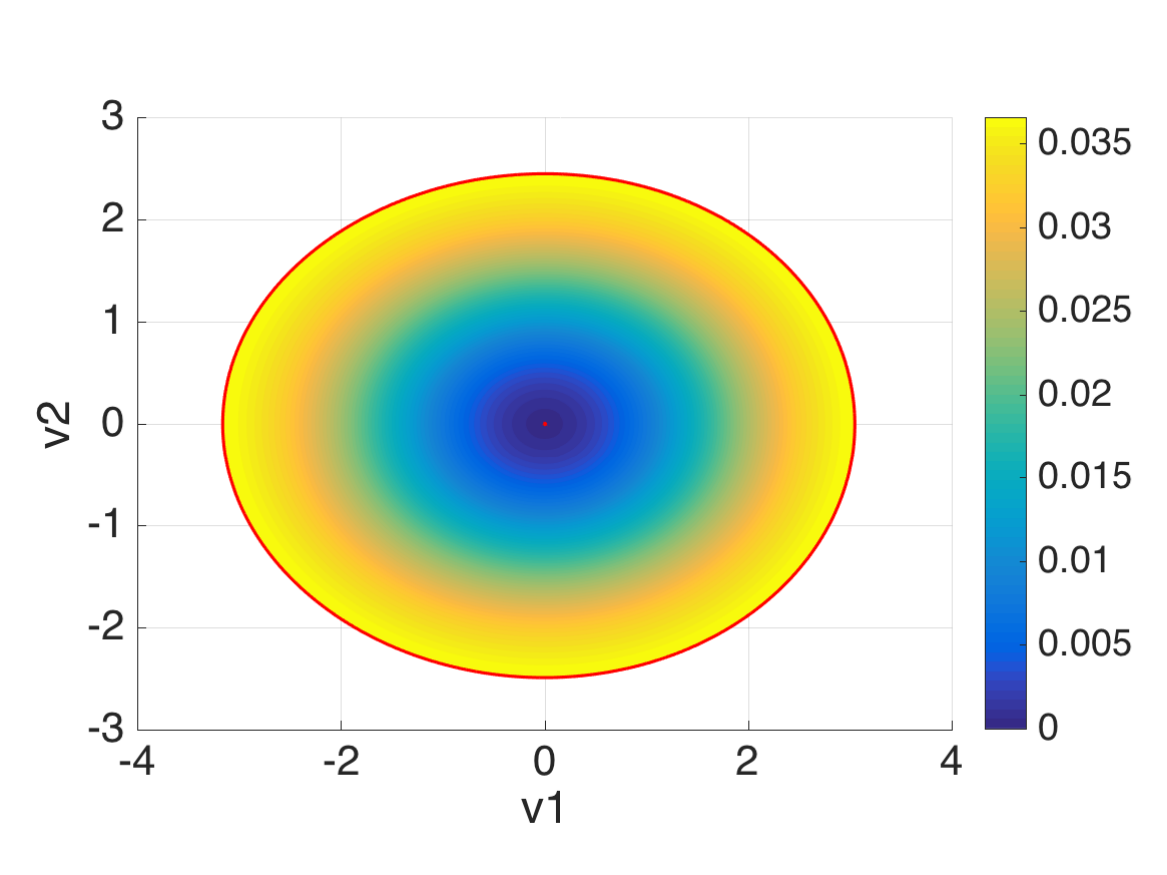}
\caption{(a): A coarsened radial mesh on the manifold $\mathcal{M}_{+}$ at $\rho=24.4$.	
The coordinate system is associated with the directions of eigenvectors of the quasipotential matrix
for the Jacobian evaluated at $C_{+}$.
(b): The quasipotential computed on this mesh.
}
\label{fig24p4Cgamma}
\end{center}
\end{figure}

We performed a 3D computation of the quasipotential with respect to $C_{+}$ aiming at obtaining the overall picture.
The computational domain was a box centered at $C_{+}$ and embracing the strange attractor.
Note that this computation is too rough to give accurate numbers, nevertheless, it captures the geometry of the level sets.
The level sets of the computed quasipotential shown in Fig. \ref{fig24p4} agree with our expectations: the quasipotential
grows until it reaches the strange attractor, it remains nearly constant on it, and then grows fast away from it, 
mainly along the union of manifolds that extends the strange attractor.
Again, we performed a 2D computation on the manifold $\mathcal{M}_{+}$ on a radial $6001\times7200$ mesh
with $K_r = 150$ and $K_a = 500$
and found the quasipotential at $\gamma_{+}$ to be equal to $0.03466$ (see Fig. \ref{fig24p4Cgamma}). 
For comparison, the 3D computation performed in a cube with size 6 centered at $C_{+}$ on a
$1001\times1001\times1001$ mesh with $K = 20$ gave the quasipotential on $\gamma_{+}$
around  0.25 which is more than 7 times larger due to the issues illustrated in Fig. \ref{fig:challenge}. 
This shows that our reduction to 2D is very important for obtaining 
accurate quasipotential barriers.
%Recall that the ratio of magnitudes of the rotational and the potential components of the vector field for the
%quasipotential decomposition of the linearized system at $C_{+}$ is $973.4$.
%For comparison,  these ratios for $\rho =12$, $\rho =15$ , and $\rho = 20$ are 15.3,  23.4, and 59.4 respectively. 

\begin{figure}[htbp]
\begin{center}
(a)\includegraphics[width = 0.45\textwidth]{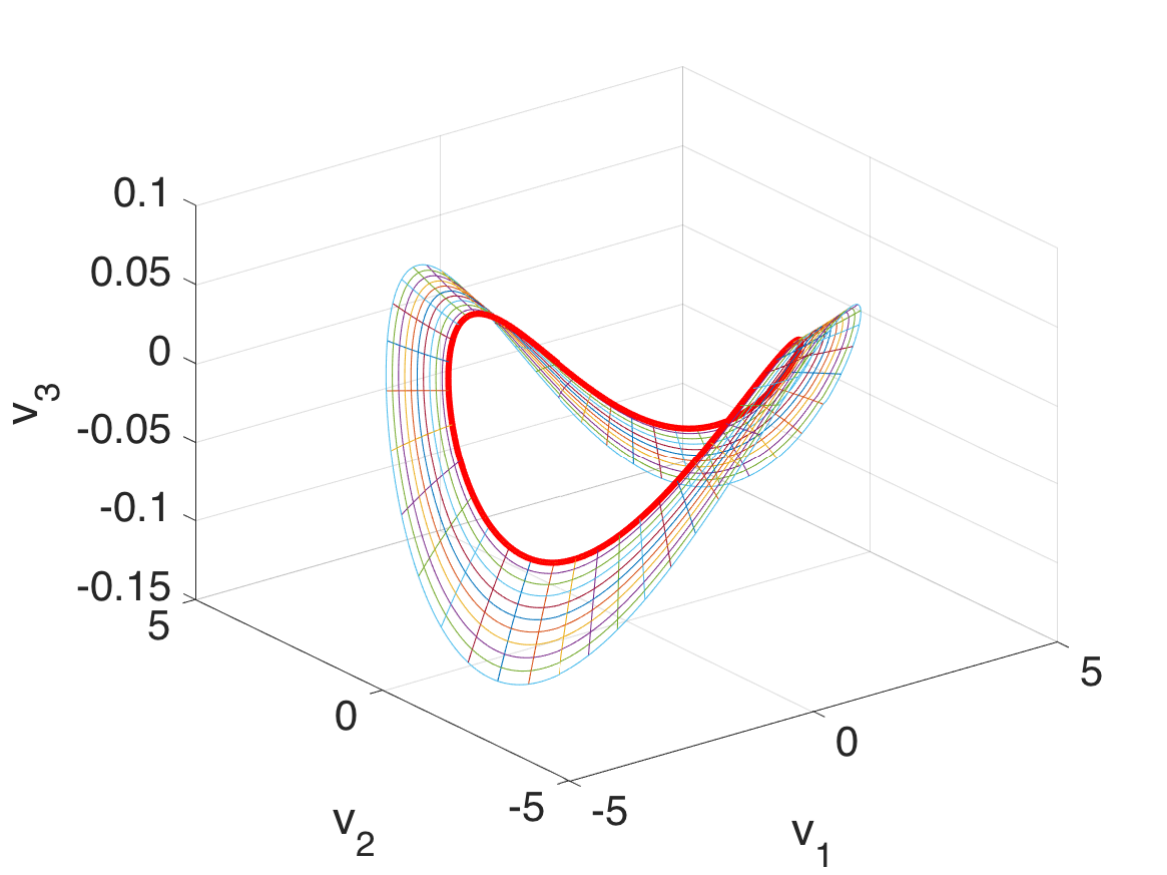}
(b)\includegraphics[width = 0.45\textwidth]{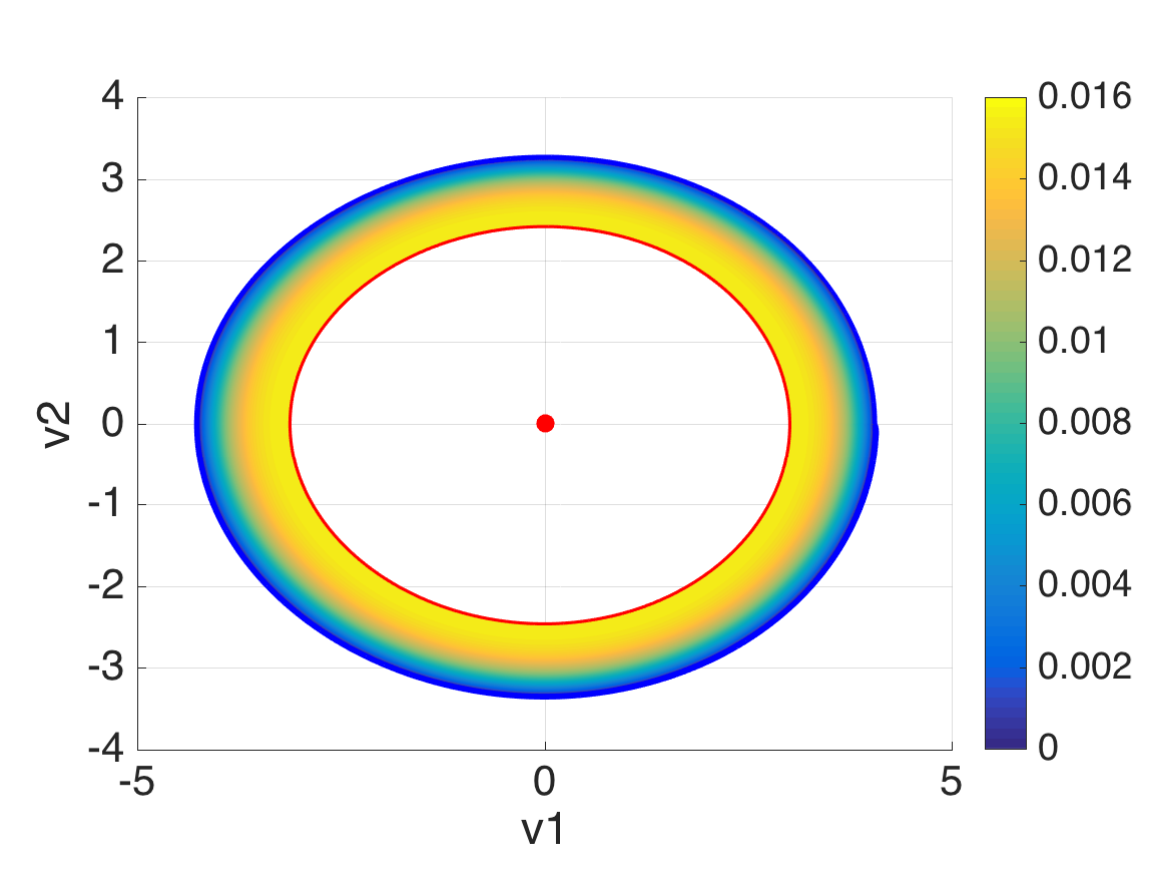}
(c)\includegraphics[width = 0.75\textwidth]{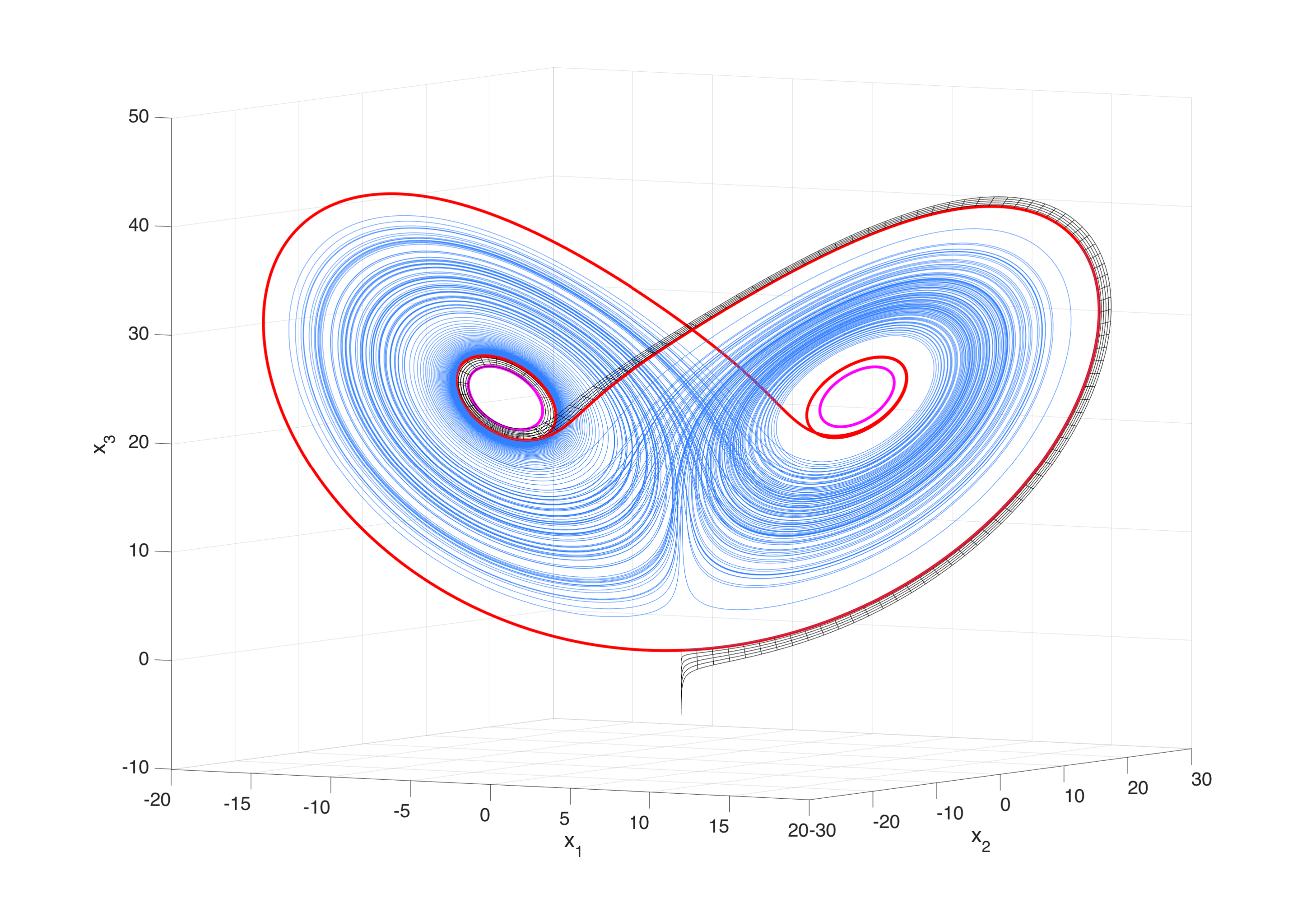}
(d)\includegraphics[width = 0.96\textwidth]{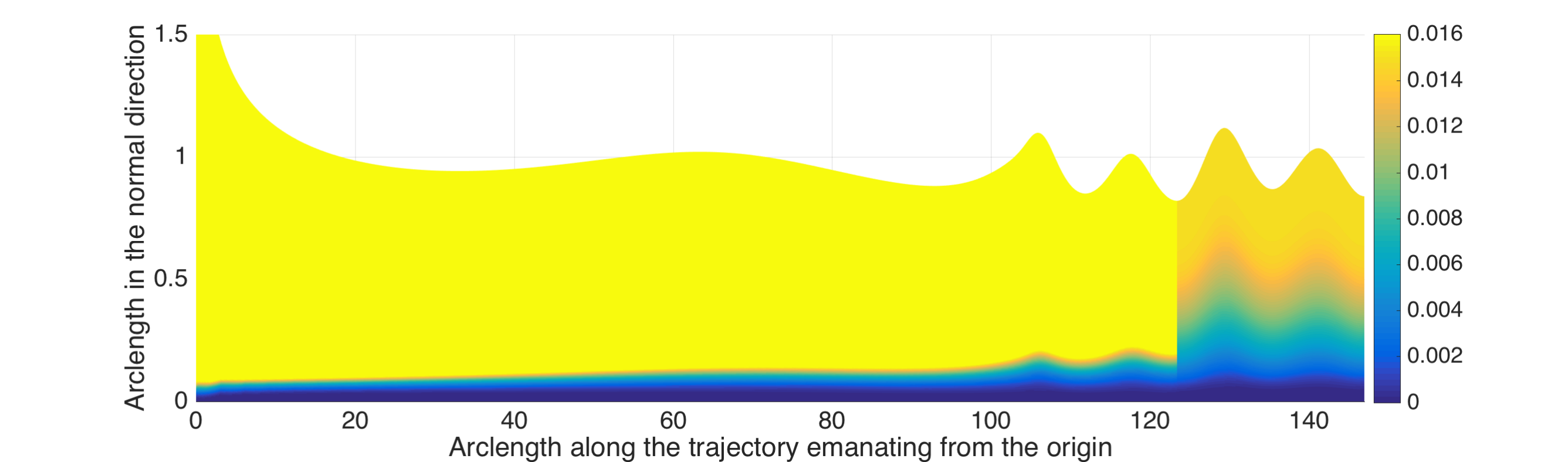}
\caption{$\rho = 24.4$. (a): A coarsened version of the``eye" mesh. The coordinate axes $\mathbf{v}_i$, $i = 1,2,3$ are 
chosen along the eigenvectors of the quasipotential matrix $Q$ of the linearized near $C_{+}$ vector field.
(b): The quasipotential computed on the ``eye" mesh.
(c): A coarsened version of the ``wing+eye" mesh. 
(d): The quasipotential computed on the ``wing+eye" mesh. 
The arclength values  less  and greater than approximately 125 correspond to the ``wing" and  ``eye" meshes respectively.
The discontinuity along the line where these meshes are glued is caused by the behavior of MAPs.
The lightest yellow region of the plot corresponds to values of the quasipotential exceeding the maximal value 0.016 on the colorbar.
}
\label{fig24p4mq}
\end{center}
\end{figure}

Fig. \ref{figSA} shows that the quasipotential level sets primarily grow along the edge of the strange attractor while remaining quite thin.
This observation suggests two possible transition mechanisms from the strange attractor to $C_{+}$.
The first one would start near the eye $Y_{+}$, climb up to $\gamma_{+}$,
and then switch to spiraling toward $C_{+}$. 
The second one would involve sliding toward $\gamma_{+}$ from the 
neighborhood of the wing $W_{-}$ to a region lying between the eye and $\gamma_{+}$
and starting spiraling toward $\gamma_{+}$ and then toward $C_{+}$. Note that a MAP for the second mechanism at $\rho=24.08$ was found in \cite{zhou}.
Coarsened versions of meshes generated for computing the quasipotential barriers for each of these transition mechanisms are
displayed in Fig. \ref{fig24p4mq} (a) and (c) respectively. 
The ``eye" mesh in Fig. \ref{fig24p4mq}(a)
is lying on the unstable loop-shaped manifold of $\gamma_{+}$ between the $\gamma_{+}$ and $Y_{+}$.
Its size is $1501\times 6000$. The found quasipotential on $\gamma_{+}$ is 0.01543 (see Fig. \ref{fig24p4mq}(b)).
The ``wing+eye" mesh  in Fig. \ref{fig24p4mq}(c)
is defined on the union of the following two manifolds.
The wing manifold is defined by trajectories starting near the negative $x_3$-semiaxis and bounded by $W_{+}$ and 
a trajectory approaching $\gamma_{+}$. The second one is the
loop-shaped unstable manifold of $\gamma_{+}$ located between $\gamma_{+}$ and $Y_{+}$.
The total mesh size is $1501\times26001$ whose $1501\times6000$ piece covers the loop.
The quasipotential computed on it is shown in Fig. \ref{fig24p4mq}(d). Its part corresponding to the loop, naturally, involves significantly smaller values 
than the one corresponding to the strip around the wing. The quasipotential value on $\gamma_{+}$ for this mesh is 0.01479 which is smaller than the one
for the eye mesh.

As we have mentioned above, the strange attractor has a finite width varying roughly from 0 to $10^{-2}$. This means that, in order to treat it as a union of four
manifolds as shown in Fig. \ref{figSA} while considering the dynamics according to SDE \eqref{sde1}, 
the parameter $\epsilon$ should be chosen at least as large as $10^{-4}$.
The discussed transition mechanisms from $A_L$ to $C_{\pm}$ are 
associated with close quasipotential barriers: the difference between them is about $5\cdot10^{-4}$.
Therefore, in order to determine which transition mechanism is dominant for $\epsilon\sim10^{-4}$, one needs to compute 
the pre-exponential factors of the  corresponding transition rates. 
Estimation of these prefactors is beyond the scope of the present work.
We leave the development of numerical methods for their evaluation for the future.

We summarize the found quasipotential barriers in Table \ref{table2}.
\begin{table}[htp]
\caption{Quasipotential barriers for stochastic Lorenz'63 \eqref{sde1} at $\sigma=10$, $\beta=\sfrac{8}{3}$, and a set of values of $\rho$.
}
\begin{center}
\begin{tabular}{|c|c|c|c|}
$\rho$ & {\bf Attractor} & {\bf Escape state} & {\bf Barrier}\\
\hline
12 & $C_{+}$ & The origin &19.5\\
15 & $C_{+}$ & $\gamma_{+}$ &18.2\\
20 & $C_{+}$ & $\gamma_{+}$ & 6.1\\
24.4 & $C_{+}$ & $\gamma_{+}$ & 0.0247\\
24.4 & $A_L$ & $\gamma_{+}$ & 0.0154 (``eye")\\
24.4 & $A_L$ & $\gamma_{+}$ & 0.0148 (``wing+eye")\\
\end{tabular}
\end{center}
\label{table2}
\end{table}%

%The second 3D computation performed in a cube with size 6 centered at $C_{+}$ on a
%$1001\times1001\times1001$ mesh with $K = 20$ gave the quasipotential on $\gamma_{+}$
%varying from 0.251 to 0.254 and averaging at 0.253.  Despite the use of large meshes,  
%our 3D computations are quite rough
%for such a complex system as Lorenz with $\rho = 24.4$. 

%%%%%%%%%%
\subsection{Perspectives and challenges for large $\rho$}
Our numerical experiments show that the level sets of the quasipotential thin out  and the diameter of the
strange attractor increases as $\rho$ grows (Fig. \ref{fig:diag}).
On one hand, this creates an underresolution problem for 3D computations
as mesh planes cannot be aligned with the level sets of the quasipotential because they are not flat.
Handling this issue by means of mesh refinement 
is limited by computer's memory.
For example, for $\rho = 100.75$ where two attracting limit cycles exist, 
the minimal level set of the quasipotential computed with respect to one of these cycles
and enclosing the other one is
thinner than the mesh step at some places.

On the other hand, thinning out of the level sets allows us to use 2D computations provided that we have 
an insight about possible transition mechanisms as
we have had for $\rho = 24.4$. This insight for larger values of $\rho$ can be gained 
from a 3D computation conducted not in a box but on a specially designed mesh.
%
%As one can see from Fig. \ref{fig:diag}, at certain intervals of values of $\rho > 24.74$, Lorenz'63 admits pairs of  attracting limit cycles. 
%For example, two symmetric limit cycles exist at $\rho = 100.75$. A question of interest 
%would be the one of finding the quasipotential barrier between them. 
%We leave this project for the future.
% because
% it would require the further enhancement of 3D quasipotential solvers
%as the level sets of the quasipotential at $\rho = 100.75$ are too thin for the current version.
%%

%%%%%%%%%%%%%%%%%%%%%%
\section{Conclusions}
\label{sec:conclusions}
We have developed a methodology for computing the quasipotential and finding quasipotential barriers 
for highly dissipative and possibly chaotic
3D dynamical systems perturbed by small white noise.
The proposed approach combines 3D computations on regular rectangular meshes with, if relevant,
dimensional reduction techniques and 2D computations on radial meshes.
%If 3D computations give accurate enough results, we are done.
%Otherwise, 3D computations are used to justify the construction of 2D manifolds or unions of manifolds for 
%more accurate 2D computations on radial meshes.
This methodology has been developed on and applied to stochastic Lorenz'63 
with $\sigma = 10$, $\beta = \sfrac{8}{3}$, and a number of values of $\rho$ ranging from 0.5 to 24.4.

%We have demonstrated an application of the recently developed 3D quasipotential 
%solver based on the ordered line integral method (OLIM) with the midpoint quadrature rule
%to the Lorenz'63 system perturbed by small white noise \eqref{sde1}. Our 3D computations have been done
%for the parameter values $\sigma = 10$, $\beta = \sfrac{8}{3}$, and a number of values of $\rho$ ranging from 0.5 to 24.4.
We have shown that, as $\rho$ increases, the level sets of the quasipotential thin out
 and the ratio of magnitudes of the rotational and potential components grows dramatically.
 On one hand, these facts render the numbers produced by 3D computations progressively less accurate.
On the other hand, the manifolds consisting of characteristics going from escape states to attractors
and those consisting of MAPs running  the other way around become very close to each other. 
This observation motivated us to
approximate the manifolds formed by the MAPs with those consisting of the characteristics. 

We have developed a technique for generating radial meshes on manifolds consisting of such characteristics
and tested our 2D OLIM quasipotential solver on an ad hoc system where the magnitude of the rotational component exceeds
that of the potential one by the factor at least as large as $10^3$, approximately as it is for $\rho=24.4$ in \eqref{sde1}.
The least squares fit for this example has given a superquadratic convergence and small 
normalized maximal absolute errors on practical mesh sizes.

Using a combination of 3D and 2D computations, we found quasipotential barriers for the escapes from the basins of
$C_{\pm}$ at $\rho = 12$, 15, 20, and 24.4. 
Furthermore, we estimated quasipotential barriers for the escape  from the basin of the Lorenz attractor at $\rho=24.4$
via two escape mechanisms. 
These barriers for 24.4 are close to each other: the difference between them is of the same order of magnitude as
the minimal  value of $\epsilon$ that makes traversing
 between different sheets of the Lorenz attractor easy.
Therefore, estimates for the pre-exponential factors for these escape rates are necessary in order to determine which transition mechanism is dominant.
We have left the development of techniques for computing these prefactors for the future.

An important advantage of computing the quasipotential in 3D is that it allows us to visualize the stochastic dynamics.
Plots of quasipotential level sets reveal the hierarchy of regions of the phase space reachable by the system perturbed by
small white noise on different timescales.  In particular, the visualization of the level sets of the
quasipotential at $\rho=24.4$ suggested us to consider and compare two possible transition mechanisms
between the strange attractor and the stable equilibria.

Our C and Matlab programs developed for the application to Lorenz'63 are posted on M. Cameron's web site \cite{mariakc} (see the package {\tt Qpot4Lorenz63.zip}) 
and on GitHub \cite{github}.

The numerical techniques developed in this work can be used for the quasipotential analysis of 
certain classes of other 2D and 3D SDEs.
The dimensional reduction to 2D can be beneficial for any 3D SDEs 
where the quasipotential with respect to an attractor grows primarily along some 2D manifold.
The use of radial meshes can dramatically improve the accuracy of 
found quasipotential thresholds in the case if the attractor 
is a stable spiral point and, perhaps, the transition state is an unstable limit cycle.

The application to the Lorenz'63 model allows us to see the limitations of the 3D quasipotential solver: 
the growth of required computational domains together with thinning out of the level sets results 
in underresolving  the latter even with the use of $1001^3$ mesh sizes. 
This motivates the directions of the future research associated with $(i)$ combining the 3D OLIMs with techniques 
for generating a 3D mesh adapted for the geometry of the problem  and 
$(ii)$ advancing the techniques for learning 2D manifolds near which the stochastic
dynamics are effectively focused.

%
%Some conclusions here.
%
%%%%%%%%%%%%%%%%%%%%%%
\section*{Acknowledgements}
We thank Dr. E. S. Kurkina for inspiring discussions on the Lorenz'63 system and sharing
numerical techniques for plotting bifurcational diagrams and
finding unstable limit cycles. We are also grateful to 
Prof. James Yorke and Prof. Kevin Lin for valuable advice regarding preparation of this manuscript.
This work is partially supported by NSF grant DMS1554907.

%
%%%%%%%%%%%%%%%%%%%%%%
%%%%%%%%%%%%%%%%%%%%%%
%
%.            A P P E N D I C E S
%
%%%%%%%%%%%%%%%%%%%%%%
%%%%%%%%%%%%%%%%%%%%%%
\appendix

%%%%%%%%%%%%%%%%%%%%%%
\section{Derivation of some equations in Section \ref{sec:background}.}
\label{sec:newapp}
\subsection{The geometric action \eqref{GA}}
Let $\phi:[T_0,T_1]\rightarrow\mr^d$ be a path with the endpoints $\phi(T_0)\in A$ and $\phi(T_1)=\mx$.
Expanding the squared norm in Eq. \eqref{FWA} and using the inequality 
$$
\|\dot{\phi}\|^2 + \|\mb(\phi)\|^2 \ge 2\|\dot{\phi}\|\|\mb(\phi)\|
$$
we obtain
\begin{equation}
\label{n1}
S_{T_0,T_1}(\phi) \ge \int_{T_0}^{T_1}\left(\|\mb(\phi)\|\|\dot{\phi} \|- \mb(\phi)\cdot\dot{\phi}\right)dt.
\end{equation}
The equality takes place if and only if 
$\|\mb(\phi)\| = \|\dot{\phi}\|$.
Since we are taking the infimum of $S_{T_0,T_1}(\phi)$ in particular with respect to $T_0$ and $T_1$, we choose 
the parametrization of $\phi$ so that $\|\mb(\phi)\| = \|\dot{\phi}\|$ and change $T_0$ and $T_1$ accordingly.
Note that $T_0$ and $T_1$ are allowed to be $-\infty$ and $+\infty$ respectively.
Next, we observe that the integral in the right-hand side of Eq. \eqref{n1}
is invariant under reparametrization of the path $\phi$. We denote the path $\phi$ reparametrized by its acrlength by $\psi$ and 
obtain Eq. \eqref{GA}.

%%%%%%%%%%%%
\subsection{The Hamilton-Jacobi equation \eqref{HJ} for the quasipotential and  equation \eqref{psi} for the MAP}
Let the path $\psi$ parametrized according to its arclength (i.e., $\|\psi'\| = 1$) be the minimizer of the geometric action \eqref{GA} 
among all absolutely continuous 
paths with one endpoint at $\mx$ and the other one at $A$.
Let us pick a small number $\delta >0$. 
Using Bellman's optimality principle \cite{bellman} and  Taylor expansion of $U$, we obtain
\begin{align*}
U(\mx)=\inf_{\|\psi'\|=1}\left\{ \int_{0}^{\delta} \left( \|\mb (\psi)\| -\mb(\psi)\cdot \psi' \right) ds +
U\left( \mx-\int_{0}^{\delta} \psi' ds\right)      \right\}\\
= \inf_{\|\psi'\|=1}\left\{ \delta \left( \|\mb (\psi)\| -\mb(\psi)\cdot \psi' - \nabla U(\mx)\cdot \psi' \right)  +
U(\mx) +O(\delta^2) \right\}.
\end{align*}
Canceling $U(\mx)$ on both sides and dividing by $\delta$ we get
\begin{align*}
0=\inf_{\|\psi'\|=1}\left\{  \|\mb (\psi)\| -\mb(\psi)\cdot \psi' -\nabla U(\mx)\cdot \psi' +  O(\delta)    \right\}.
\end{align*}
Taking the limit as $\delta \rightarrow 0,$ we obtain
\begin{align}
\label{HJder}
\inf_{\|\psi'\|=1}\left\{  \|\mb (\mx)\| - \left(\mb(\mx)+   \nabla U(\mx) \right)\cdot \psi'  \right\}=0.
\end{align}
The infimum is attained when the term
$\left(\mb(\mx)+  \nabla U(\mx) \right)\cdot \psi' $
is maximal, i.e., when
\begin{equation}
\label{psimax}
 \psi'=\displaystyle\frac{\mb(\mx) + \nabla U(\mx)}{\| \mb(\mx) + \nabla U(\mx) \|}.
\end{equation}
Observing that $\mx$ is the point of the path $\psi$ at which $\psi'$ is evaluated, we see that \eqref{psimax} coincides with equation \eqref{psi}.
Plugging  \eqref{psimax} into \eqref{HJder}, we get
\begin{equation}
\label{HJder1}
\|\mb(\mx)\|  = \| \mb(\mx) + \nabla U(\mx) \|.
\end{equation}
Taking squares of both sides of Eq. \eqref{HJder}, canceling $\|\mb(\mx)\|^2$, and dividing by 2, we obtain the 
desired Hamilton-Jacobi equation \eqref{HJ}:
\begin{equation*}
\frac{1}{2}\|\nabla U(\mx)\|^2 +  \mb(\mx)\cdot \nabla U(\mx)=0.
\end{equation*}

%%%%%%%%%%%%%%%%%%%%%%%%

\section{The dynamics of the Lorenz system \eqref{lorenz}}
\label{appLorenz}

Let us fix the parameters $\sigma = 10$ and $\beta=\sfrac{8}{3}$.  
As $\rho$ grows from zero to infinity, the dynamics of \eqref{lorenz} go through a number of bifurcations \cite{preturbulence,sparrow,sparrowIEEE,strogatz}.
\begin{figure}[htbp]
\begin{center}
\includegraphics[width = 0.85\textwidth]{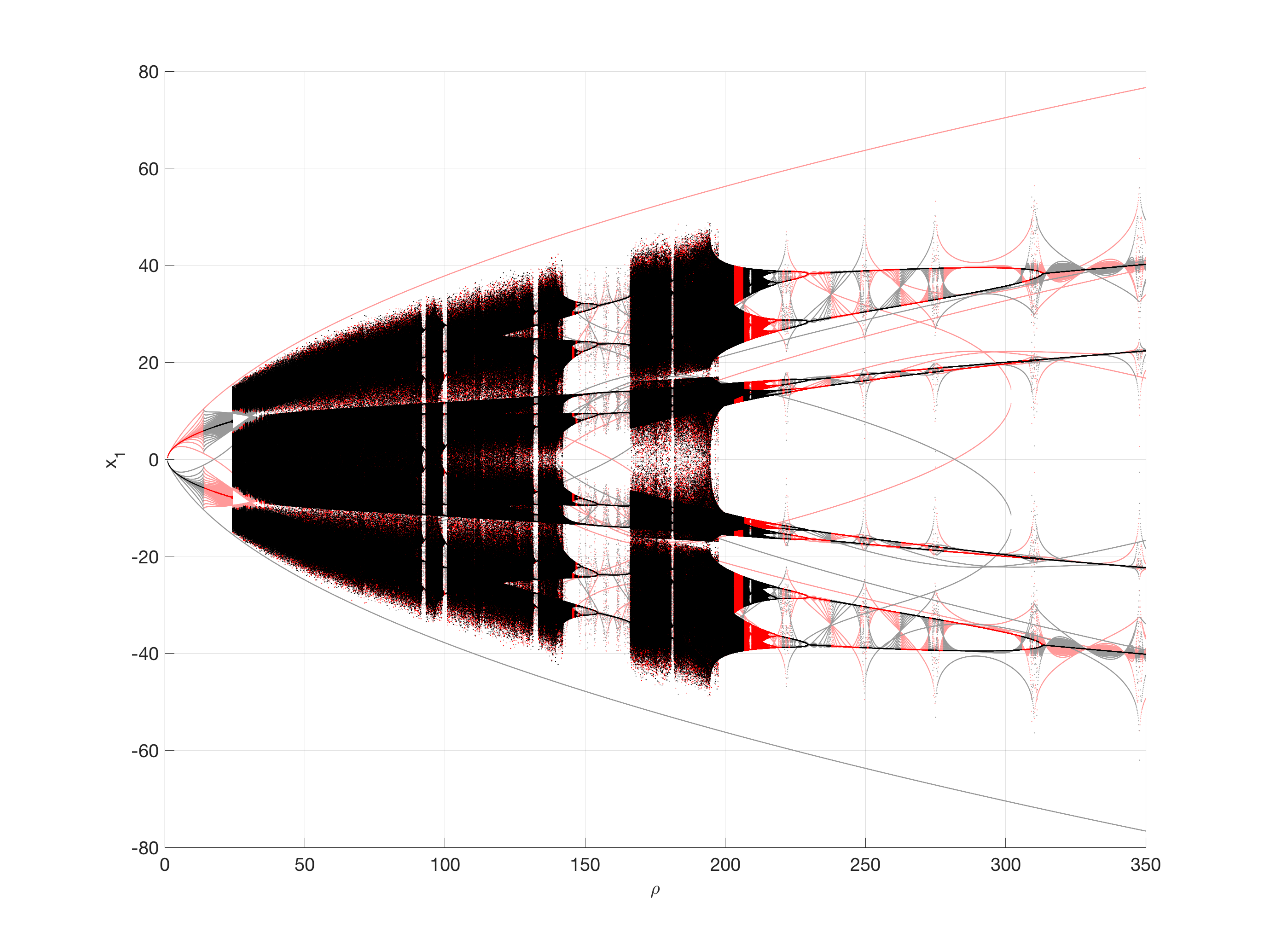}
%(b)\includegraphics[width = 0.85\textwidth]{diagram_zoom.pdf}
\caption{Consider the characteristics of \eqref{lorenz} emanating from the origin along the directions $\mathbf{\xi}$ and $-\mathbf{\xi}$ and traced
for the time interval $0\le t\le 200$.
The $x_1$-components of their intersections with the horizontal plane passing through the equilibria $C_{\pm}$ are plotted for $1\le \rho\le 350$
with pink and grey dots respectively. Then each characteristic is continued to be traced for $200\le t\le 400$.
The resulted $x_1$ components of their intersections with the same plane are marked with red and black respectively.
The dashed green vertical lines correspond to 
the critical values of $\rho$: $\rho_0\approx 13.926$, $\rho_1\approx 24.06$, and $\rho_2\approx 24.74$.
}
\label{fig:diag}
\end{center}
\end{figure}
\begin{figure}[htbp]
\begin{center}
\includegraphics[width = 0.8\textwidth]{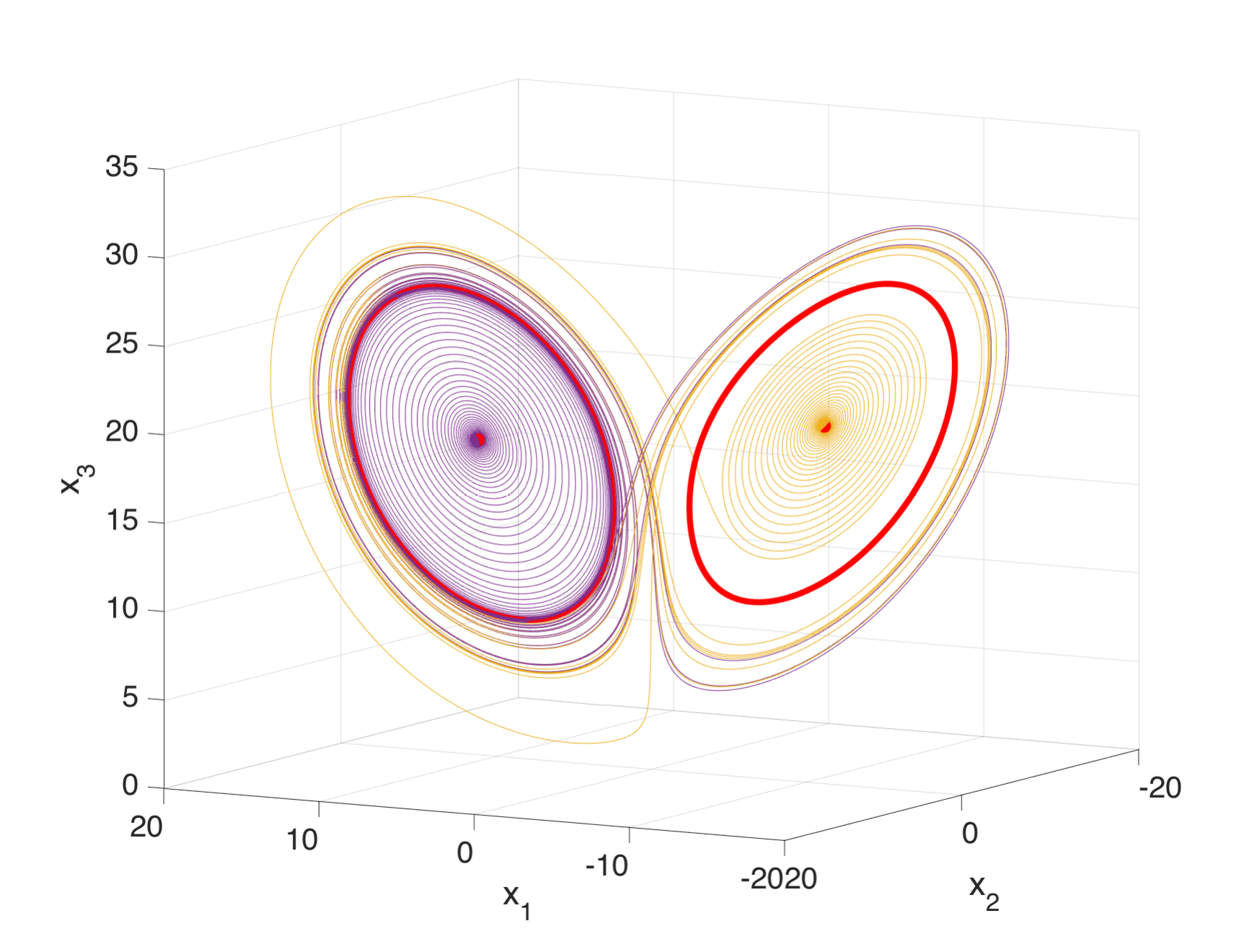}
\caption{(a): An example of two characteristics at $\rho = 20$ starting at two close points lying near $\gamma_{+}$ 
on the cone with vertex at $C_{+}$ and consisting 
of all rays passing through $\gamma_{+}$ and eventually diverging and approaching different equilibria.
}
\label{fig:pret}
\end{center}
\end{figure}

\begin{itemize}
\item For all $0<\rho <\infty$, the origin is a fixed point of \eqref{lorenz}.
It is the only equilibrium for $0<\rho < 1$, and it is globally attracting. 
At $\rho =1$, a supercritical pitchfork bifurcation occurs transforming the origin into a Morse index one saddle and 
giving birth to two equilibria
\begin{equation}
\label{Cpm}
C_{\pm} = \left(\pm\sqrt{\beta(\rho-1)},\pm\sqrt{\beta(\rho-1)},\rho-1\right).
\end{equation}
They remain asymptotically stable for $1<\rho<\rho_2\approx 24.74$.
The unstable manifold of  \eqref{lorenz} linearized near the saddle at the origin for $1 < \rho < \infty$
is the span of the vector
\begin{equation}
\label{uvec}
\mathbf{\xi}=\left[\begin{array}{c}\sigma\\\frac{\sigma - 1}{2} +\sqrt{\left(\frac{\sigma + 1}{2}\right)^2 +\sigma(\rho - 1)}\\0\end{array}\right].
\end{equation} 
To delineate the evolution of the dynamics of \eqref{lorenz} as $\rho$ grows from 1 to infinity, we have plotted the bifurcation diagram
displayed in Fig. \ref{fig:diag}. For each $\rho$ from  $1.05$ to $349.95$ with step $0.1$, we traced the trajectory starting at $10^{-2}\xi$ for time $0\le t\le 200$
and recorded its points of intersection with the plane 
$$
\alpha = \{\mx~|~x_3 = \rho - 1\}
$$ 
passing through the equilibria $C_{\pm}$. The $x_1$-components of these intersects
are shown with pink dots in the $(\rho,x_1)$-plane. The time interval $0\le t\le 200$ is large enough for this trajectory to approach an attractor. 
Then, in order to depict $x_1$-components of the intersection of the attractor with the plane $\alpha$, we continued tracing the trajectory for $200\le t\le 400$
and plotted the $x_1$ components of its intersects with $\alpha$ with red dots. The corresponding sets of points for the trajectory starting at $-10^{-2}\xi$
are obtained using the aforementioned symmetry of \eqref{lorenz}. They are plotted with grey and black dots respectively. This procedure is implemented 
in the Matlab code \verb|lorenz_diagram.m|.

\item For $1 < \rho < \rho_0\approx 13.926$, 
the characteristics emanating from the saddle at the origin along the directions $\mathbf{\xi}$ and $-\mathbf{\xi}$
approach, respectively, $C_{+}$ and $C_{-}$ without crossing the plane $x_1 = 0$ (see Fig. \ref{fig:diag}).

\item The interval $13.926 \approx \rho_0  < \rho < \rho_2\approx 24.74$
is marked by the existence of the saddle limit cycles $\gamma_{+}$ and $\gamma_{-}$ surrounding $C_{+}$ and $C_{-}$ respectively.
The equilibria $C_{\pm}$ remain the only attractors for $\rho_0 < \rho <\rho_1\approx 24.06$.
At $\rho=\rho_0$, there exist homoclinic orbits 
emanating from the origin and approaching it as $t\rightarrow\infty$.
For all $\rho_0 <\rho < \rho_1$, the characteristics emanating from the origin along the directions $\mathbf{\xi}$ and $-\mathbf{\xi}$
go approximately half-way around the limit cycles, cross the plane $x_1=0$, and approach $C_{-}$ and $C_{+}$ respectively (see Fig. \ref{fig:diag}).
As $\rho$ grows within this interval,  there  develops a phenomenon called \emph{preturbulence} \cite{preturbulence},
characterized by chaotic behavior and divergence of close characteristics in a region surrounding $\gamma_{\pm}$. 
Let $\Upsilon_{+}$ be a cone consisting of all rays
starting at $C_{+}$ and crossing $\gamma_{+}$, i.e.,
\begin{equation}
\label{cone}
\Upsilon_{+}: = \{ C_{+} + t(\mx - C_{+}) ~|~t\ge 0,~\mx\in\gamma_{+} \}.
\end{equation}
Characteristics starting  on $\Upsilon_{+}$ near and outside $\gamma_{+}$ 
perform more and more revolutions around $C_{+}$ and $C_{-}$ prior they settle to 
spiraling near one of the stable equilibria. Moreover,  as $\rho$ tends to $\rho_1$,
it is getting progressively harder and finally impossible to predict using double-precision arithmetic 
which equilibrium such a characteristic 
will eventually approach.
An example of two characteristics for $\rho = 20$ starting at two close points near $\gamma_{+}$ on the  cone $\Upsilon_{+}$ 
and eventually approaching different equilibria  is shown in Fig. \ref{fig:pret}.
At $\rho = \rho_1$, the characteristics emanating from the origin along the directions $\mathbf{\xi}$ and $-\mathbf{\xi}$
approach $\gamma_{-}$ and $\gamma_{+}$ respectively. This gives birth to a strange attractor a.k.a. the Lorenz attractor.
We will denote it by $A_L$.

% where ``$L$" stands for ``Lorenz".
%The critical value $\rho_1\approx24.06$ was found in \cite{yorke} by means of fitting
%a certain class of functions to the Lorenz map. 
%
\item For $ 24.06\approx\rho_1<\rho<\rho_2\approx24.74$, 
there are three attractors: the strange attractor $A_L$, and the asymptotically stable equilibria $C_{\pm}$.
The characteristics  emanating from the origin along $\pm\mathbf{\xi}$ 
miss the saddle cycles $\gamma_{\mp}$ respectively and start spiraling away from them.
$\gamma_{\pm}$ lie on the boundaries of the basins of $C_{\pm}$ respectively
and, as we show in Section \ref{sec:coex}, play roles of the escape states. %  from $C_{\pm}$ to $A_L$.
At $\rho = \rho_2$, the saddle cycles $\gamma_{\pm}$ shrink to the corresponding equilibria  $C_{\pm}$
rendering them unstable, i.e., a subcritical Hopf bifurcation takes place.
\item For $ 24.74\approx\rho_2<\rho < \infty$ , the dynamics are complicated as can be inferred from Fig. \ref{fig:diag}. 
$A_L$ is the only attractor for some open interval of $\rho$ starting at $\rho_2$ (Fig. \ref{fig:diag}). 
It exists for a union of intervals of $\rho$ stretching up to approximately $\rho =215.364$ \cite{sparrow}.
The interval $\rho_2<\rho\lesssim 215.364$ is cut through by a number of windows of periodicity where there exist attracting limit cycles.
The largest of them is $145\lesssim\rho\lesssim166$. Other windows are seen around $\rho =93$, $\rho = 100$, $\rho = 133$,
and $\rho = 181.5$.
Zooming in, we can spot more windows of periodicity (see Fig. \ref{fig:diag}) and reveal
cascades of period doublings marking the Feigenbaum scenarios of transition to chaos.
The final doubling period interval  $215.364\lesssim \rho\lesssim 313$  \cite{sparrow} is clearly visible in Fig. \ref{fig:diag}. 
Near $\rho=313$, two symmetric attracting limit cycles cycles
merge into one resulting in the final limit cycle that remains the only attractor for all larger values of $\rho$. 
\end{itemize}

%We will quantify and visualize the dynamics of SDE \eqref{sde1} 
%by means of computing the quasipotential for a number of representative values of $\rho$.

%%%%%%%%%%%%%%%%%%%%%%

\section{The KKT conditions for the simplex update}
\label{appKKT}
The Lagrange function for the constrained minimization problem \eqref{3ptu}--\eqref{3con} is
\begin{equation}
\label{lagf}
\mathcal{L}(\lambda,\mu) = U_{\lambda} + \mathcal{Q}_M(\mx_{\lambda},\mx) 
-\mu_1\lambda_1-\mu_2\lambda_2-\mu_3(1-\lambda_1-\lambda_2),
\end{equation}
where $\lambda = [\lambda_1,\lambda_2]$ and $\mu = [\mu_1,\mu_2,\mu_3]$.
For brevity, we denote the function to be minimized by $f$:
$$
f(\lambda):= U_{\lambda} + \mathcal{Q}_M(\mx_{\lambda},\mx).
$$
The KKT optimality conditions applied to \eqref{lagf} are
\begin{align}
\nabla_{\lambda}\mathcal{L}(\lambda,\mu) & = \nabla f(\lambda)
-\mu_1\left[\begin{array}{c}1\\0\end{array}\right] 
-\mu_2\left[\begin{array}{c}0\\1\end{array}\right] 
-\mu_3\left[\begin{array}{c}-1\\-1\end{array}\right] = \left[\begin{array}{c}0\\0\end{array}\right] \label{kkt1}\\
\mu_1 &\ge 0,\quad \mu_2\ge 0,\quad \mu_3\ge 0,\label{kkt2}\\
\lambda_1&\ge 0,\quad \lambda_2\ge 0,\quad 1-\lambda_1-\lambda_2\ge 0,\label{kkt3}\\
\lambda_1\mu_1 &= 0,\quad \lambda_2\mu_2 = 0,\quad (1-\lambda_1-\lambda_2)\mu_3 = 0. \label{kkt4}
\end{align}
Let us check whether the initial guess $\lambda=[\lambda^{\ast},0]$ where 
$\lambda^{\ast}$ is the minimizer of $f$ on $[\lambda_1,0]$, $0<\lambda_1<1$, 
corresponding to the line segment $[\mx_0,\mx_1]$, satisfies the KKT conditions \eqref{kkt1}--\eqref{kkt4}.
Condition \eqref{kkt4} with $\lambda_1=\lambda^{\ast}\in(0,1)$ and $\lambda_2=0$ implies that $\mu_1=\mu_3=0$.
Therefore, the first component in \eqref{kkt1} is zero as 
\begin{equation}
\label{f1ast}
\frac{\partial}{\partial\lambda_1}f(\lambda^{\ast},0) = 0.
\end{equation}
The second component of \eqref{kkt1} must be also zero, hence 
\begin{equation}
\label{f2ast}
\frac{\partial}{\partial\lambda_2}f(\lambda^{\ast},0) -\mu_2 = 0.
\end{equation}
Condition \eqref{kkt2} demands that $\mu_2\ge 0$. Hence, $\lambda=[\lambda^{\ast},0]$ is a solution of 
the constrained minimization problem \eqref{3ptu}--\eqref{3con} if 
\begin{equation}
\label{f2astast}
\mu_2 = \frac{\partial}{\partial\lambda_2}f(\lambda^{\ast},0) \ge 0,
\end{equation}
i.e., if equation \eqref{kktcheck} holds. In this case, we reject the simplex update. 
Otherwise, we proceed with solving the minimization problem
 \eqref{3ptu}--\eqref{3con}.

%%%%%%%%%%%%%%%%%%%%%%

\section{Quasipotential decomposition for linear SDEs}
\label{sec:app0}
In this Appendix, we explain how one can find the quasipotential for 
linear SDEs for which the origin is an asymptotically stable equilibrium.
This is useful for initializing the OLIMs near asymptotically stable equilibria and
for estimating the ratio of the magnitudes of the rotational and potential components of the vector field.

Let $J$ be a $d\times d $ matrix with all eigenvalues having negative real parts. 
In this work, $J$ is the Jacobian matrix
of the vector field $\mb$ evaluated at an asymptotically stable equilibrium $\mx^{\ast}$ of $\dot{\mx} = \mb(\mx)$.
We consider the linear SDE for the variable $\my: = \mx - \mx^{\ast}$:
\begin{equation}
\label{linSDE}
d\my = J\my dt + \sqrt{\epsilon}d\mw.
\end{equation}
The problem of finding the quasipotential decomposition for the vector field $J\my$ reduces to the problem of finding a
symmetric positive definite matrix $Q$ such that \cite{chen,chen1}
\begin{equation}
\label{lord}
\my^\top Q (J + Q)\my = 0~~\text{ for all}~~ \my\in\mr^d.
\end{equation} 
The matrices $Q$ and  $L: = J + Q$ are  called the \emph{quasipotential matrix} and the \emph{rotational matrix} respectively.
Condition \eqref{lord} is equivalent to the requirement that the matrix $Q(J + Q)$ is antisymmetric,
i.e., $Q(J + Q) + (J + Q)^\top Q = 0$. The last equation for $Q$ is reducible to a Sylvester equation for $Q^{-1}$ 
and has a unique positive definite solution
that can be found using the Bartels-Stewart algorithm implemented in Matlab in the command {\tt sylvester} (see \cite{YPC} for details).

To make our quasipotential solver for the Lorenz system self-contained and 
facilitate experiments with various values of $\rho$, 
we have developed a C code {\tt LinLorenz.c} for finding the quasipotential 
decomposition for the Lorenz system linearized near its asymptotically stable equilibria.
The quasipotential decomposition is found by an algorithm similar to Bartels-Stewart but simplified and customized for Lorenz'63.
A description of it is linked to the provided software package \cite{mariakc}.

Once the quasipotential decomposition for a linearized system is available, one can obtain an estimate for the ratio $\Xi(\mx)$
of the magnitudes of the rotational and potential components near asymptotically stable equilibria:
\begin{equation}
\label{xiest}
\Xi\lesssim\max_{\|\my \|= 1}\frac{\|L\my\|}{\|Q\my\|}.
\end{equation}
The graph of the right-hand side of \eqref{xiest} with $J$ been the Jacobian matrix evaluated at $C_{+}$ of \eqref{lorenz}
is plotted in Fig. \ref{fig:rotpot} for the range $1<\rho<\rho_2\approx 24.74$.

%%%%%%%%%%%%%%%%%%

\section{Proof of Lemma \ref{lemma1}}
\label{appL1}
\begin{proof}
First we prove that the manifold $\mathcal{M}'$ consisting of MAPs going from the attractor $A$ to the curve $\gamma$ lies 
in the sublevel set $\mathcal{V}_a$. 
Let $\psi$ be a MAP going from $A$ to $\gamma$. 
Since $\mathcal{V}_a$ completely lies in the basin of $A$, the quasipotential strictly increases along the MAP. 
Therefore, for any $\my$ lying on the path  $\psi$,
$U(\my) \le a$ which means that $\psi\subset\mathcal{V}_a$. Since this is true for all such MAPs, $\mathcal{M}'\subset\mathcal{V}_a$.

Now let us prove that the manifold $\mathcal{M}$ consisting of all characteristics starting at $\gamma$ and running to $A$ 
lies in $\mathcal{V}_a$. We proceed from converse. Suppose a characteristic starting at $\gamma$ and going to $A$ leaves 
$\mathcal{V}_a$ at a point $\mx_0$ and reenters $\mathcal{V}_a$ at a point $\mx_1$ after that. Let $\my$ be a point of this characteristic 
located between $\mx_0$ and $\mx_1$. Since the motion of the characteristic contributes nothing to the Freidlin-Wentzell action \eqref{FWA},
 $U(\my) = U(\mx_0) = a$. This contradicts to the assumption that $\my\notin\mathcal{V}_a$. Therefore, the characteristic must completely lie in
 $\mathcal{V}_a$. Since this argument applies all characteristics constituitng $\mathcal{M}$, 
 we conclude that $\mathcal{M}\subset\mathcal{V}_a$.
\end{proof}

%%%%%%%%%%%%%%%%%
%
\section{Building radial meshes}
\label{sec:appA}
Suppose we would like to build a radial mesh on a 2D manifold formed 
by characteristics of $\dot{\mx}=\mb(\mx)$
going from an unstable limit cycle $\gamma$ to an asymptotically stable spiral point $\mx^{\ast}$.
First, we pick a set of points $\mx^k$, $k = 0,1,\ldots,N_a-1$, equispaced along $\gamma$.
For each point $\mx^k$, we define a plane $\alpha^k$ passing through
$\mx^{\ast}$ and  $\mx^k$ whose normal $\mathbf{a}^k$ lies in the plane spanned by $\mb(\mx^k)$ and $\mx^k-\mx^{\ast}$.

Then,  we trace
a trajectory $\my(t)$ starting near $\gamma$ and ending upon reaching a $\delta$-ball centered at $\mx^{\ast}$ where $\delta$ is a small number.
Let $\my_1,\ldots, \my_n$ be the set of intersects of $\my(t)$ with the plane $\alpha^0$ 
at which the sign of $(\my(t) - \mx^0)^\top \mathbf{a}^0$ changes from ``$-$" to ``$+$".
Adding $\mx^0$ and $\mx^{\ast}$ to this set and interpolating, we get a curve lying in $\alpha_0$
and connecting $\gamma$ and $\mx^{\ast}$. We define a set of points $\{\mz^0_i\}_{i=0}^{N_r-1}$ uniformly  
distributed along this curve such that $\mz^0_0\equiv \mx^{\ast}$ and $\mz^0_{Nr-1}\equiv\mx^0$.

Next, for $k=0,1,2,\ldots,N_a-2$, we trace the trajectories starting at $\mz_i^k$, $i=1,\ldots,N_r-2$, and terminate them as soon as they 
reach the plane $\alpha_{k+1}$. 
As above, we add $\mx_{k+1}$ and $\mx^{\ast}$ to these terminal points, interpolate them, and 
pick a set of points $\mz_i^{k+1}$, $i = 0,\ldots,N_r-1$, uniformly distributed along the interpolant and such that 
$\mz_0^{k+1}\equiv \mx^{\ast}$ and $\mz^{k+1}_{N_r-1}\equiv \mx^{k+1}$.
As a result, we obtain the radial mesh
$$
\{\mz_i^k~|~0\le i\le N_r - 1,~0\le k\le N_a-1\}.
$$
This procedure is implemented in the Matlab code {\tt make2Dmesh.m}\\
 in  {\tt Qpot4Lorenz63.zip} \cite{mariakc}.

Similar methodologies have been used to construct radial meshes between two simple closed curves and 
between two given segments of two  distinct characteristics.

%%%%%%%%%%%%%%%%%%%%%%%
%Appendix B: thickness of the strange attractor

\section{Estimating the width of the Lorenz attractor} 
\label{sec:appB}
\begin{figure}[htbp]
\begin{center}
\includegraphics[width = 0.7\textwidth]{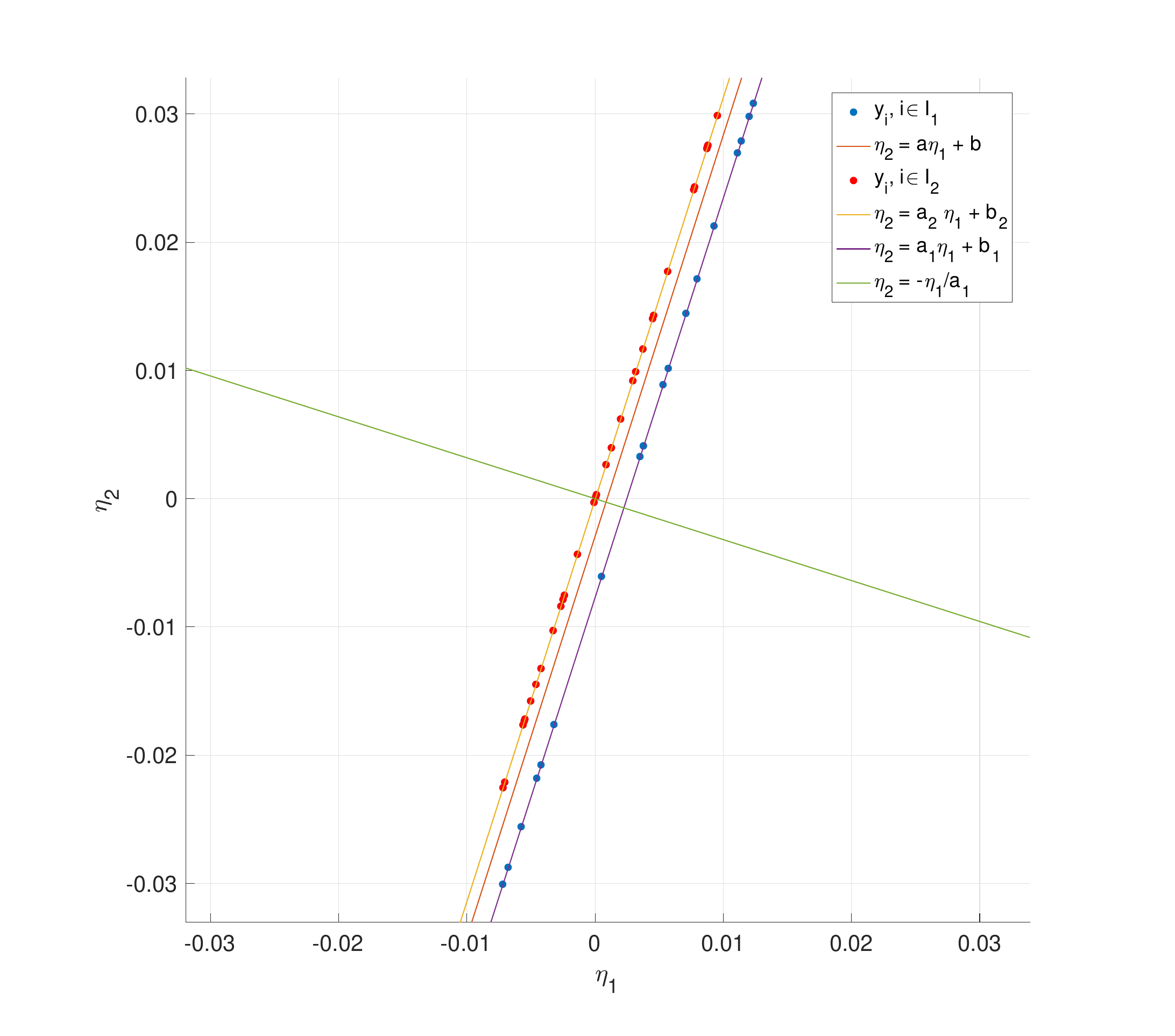}
\caption{Estimating the thickness of the Lorenz attractor using linear least squares fits in a Poincare section.}
\label{fig:thmeas}
\end{center}
\end{figure}

Let $\mx$ be a point lying on the Lorenz attractor $A_L$ and let $\alpha$ be
the plane passing through $\mx$ and normal to $\mb(\mx)$ where $\mb$ is the Lorenz vector field, i.e.,
$$
\alpha:=\{\mathbf{z}\in\mr^3~|~(\mathbf{z} - \mx)^\top\mb(\mx) = 0\}.
$$
We trace a trajectory $\my(t)$ starting at $\mx$ for time 
$10^4$ and record the points $\my_i$, $1\le i \le N$, at which the sign of $(\my(t)-\mx)^\top\mb(\mx)$
switches from ``$-$" to ``$+$".
We set up a Cartesian coordinate system $(\eta_1,\eta_2)$ in the plane $\alpha$ with the origin at $\my_1\equiv \mx$
and find the coordinates of the recorded points $\my_i$: $\my_i\equiv(\eta_1^i,\eta_2^i)$.
We pick a square $S:=[-0.25\le\eta_1\le0.25]\times[-0.25\le\eta_2\le0.25] $ in this plane and select the subset $I\subset\{1,\ldots,N\}$ 
such that the points $\my_i$, $i\in I$, lie in $S$.
Visualizing the set $\my_i$, $i\in I$, and zooming in if necessary, we see  that they are arranged near two almost parallel lines (see Fig. \ref{fig:thmeas}).
The least squares fit to this set of points with a linear function $\eta_2 = a\eta_1 + b$
gives a line dividing it into two subsets: 
\begin{align*}
I_1 &= \{i\in I~|~\eta_2^i < a\eta_1^i + b\}~~{\rm and} \\
I_2 &= \{i\in I~|~\eta_2^i > a\eta_1^i + b\}.
\end{align*}
Next, we find linear least squares fits $\eta_2=a_1\eta_1 + b_1$ and $\eta_2=a_2\eta_1 + b_2$
for the subsets of $\my_i$ corresponding to $ I_1$ and $I_2$, respectively.
One of these linear functions must pass very close to the origin because $\mx$ lies near one of these lines, hence either $b_1$ or $b_2$ 
is very close to zero in comparison with the other one. Assume that $|b_2|\ll|b_1|$. 
If this is the other way around, we swap the notations.
Also, these lines are almost parallel, hence $a_1$ and $a_2$ are very close.
Finally, we find a line orthogonal to $\eta_2=a_1\eta_1 + b_1$ and passing through the origin: $\eta_2 = -a_1^{-1}\eta_1$.
Then the thickness of $A_L$ near $\mx$ is approximately equal to the distance between the origin and 
the intersect of $\eta_2 = -a_1^{-1}\eta_1$ and
$\eta_2=a_1\eta_1 + b_1$.
This technique is implemented in the Matlab program {\tt thickness.m} 
\cite{mariakc,github}.

\bibliographystyle{siamplain}

%--- Bibliography ---%

%Journal of Experimental and Theoretical Physics
%June 2007, Volume 104, Issue 6, pp 966?971 | Cite as
%Channels and jokers in continuous systems
%Authors
%Authors and affiliations
%O. Ya. Butkovski?M. Yu. Logunov
%1.
%Statistical, Nonlinear, and Soft Matter Physics
%Received: 16 October 2006
%\bibliography{references}
\end{document}